\documentclass{article}

\usepackage{amsfonts}
\usepackage{amsmath}
\usepackage{latexsym}
\usepackage{amssymb}
\usepackage{graphics}
\usepackage{graphicx}

\usepackage{caption}
\usepackage{subcaption}

\usepackage{algorithm}
\usepackage{algpseudocode}
\usepackage{fancybox}
\usepackage{float}

\newcommand\Algphase[1]{%
	\Statex\hspace*{-\algorithmicindent}\textbf{#1}%
}

\evensidemargin=\oddsidemargin
\newtheorem{ass}{Assumption}
\newtheorem{theorem}{Theorem}[section]
\newtheorem{proposition}[theorem]{Proposition}
\newtheorem{lemma}[theorem]{Lemma}
\newtheorem{definition}[theorem]{Definition}
\def\la{\lambda}
\def\eps{\epsilon}
\def\Re{\mathbb{R}}
\def\one{\mathbf{1}}

\title{A derivative-free approach to mixed integer constrained multiobjective nonsmooth optimization}

\author{
G. Liuzzi\footnotemark[1] \and 
S. Lucidi\footnotemark[2]
}

\date{\today}
\begin{document}

\maketitle

\renewcommand{\thefootnote}{\fnsymbol{footnote}}

\footnotetext[1]{``Sapienza'' Universit\`a  di Roma, Dipartimento di Ingegneria Informatica Automatica e
Gestionale ``A. Ruberti'', Via Ariosto 25, 00185 Rome, Italy and Consiglio Nazionale delle Ricerche, Istituto di Analisi dei Sistemi e Informatica ``A. Ruberti'', Via dei Taurini 19, 0085 Rome, Italy. \texttt{liuzzi@diag.uniroma1.it}}
\footnotetext[2]{``Sapienza'' Universit\`a  di Roma, Dipartimento di Ingegneria Informatica Automatica e
Gestionale ``A. Ruberti'', Via Ariosto 25, 00185 Rome, Italy. \texttt{lucidi@dis.uniroma1.it}}

\renewcommand{\thefootnote}{\arabic{footnote}}

{\small
\begin{abstract}
In this work, we consider multiobjective optimization problems with both bound
constraints on the variables and general nonlinear constraints, where objective and constraint function
values can only be obtained by querying a black box. Furthermore, we consider the case 
where a subset of the variables can only take integer values. We propose a new linesearch-based 
solution method and show that it converges to a set of stationary points for the problem.
For what concerns the continuous variables, we employ a strategy for the estimation of the Pareto frontier
recently proposed in the literature and which takes advantage of dense seuqnces of search directions.
The subset of variables that must assume dicrete values are dealt with using primitive directions 
appropriately modified to take into account the presence of more than one objective functions.
Numerical results obtained with the proposed method on a set of test problems and comparison
with other solution methods show the viability and efficiency of the proposed approach.
\end{abstract}
}

\par\medskip

{\bf AMS}
90C30, 90C11, 90C56, 65K05, 49J52


\section{Introduction}
In many applications (especially coming from engineering or economic contexts) the simultaneous minimization of several conflicting objective functions is to be considered within a feasible region described by (possibly) nonlinear constraints. It is often the case that objective and/or constraint function values are computed by using complex simulation programs so that first (or higher) order information is not available. Furthermore, function values are typically affected by a certain degree of numerical noise which prevents the use of approximation techniques to compute derivatives. A further distinguishing feature of many multiobjective problems that arise in applications is the presence of variables that can only assume integer or discrete values. 

The aim of the paper is to develop and analyse a derivative-free method for the solution of multiobjective mixed integer constrained optimization problems. More specifically, we are interested in the definition of an algorithm belonging to the calss of ``{\em a posteriori}'' methods, namely those algorithms that try to compute an approximation of the Pareto front for the problem. Indeed, when multiple objective functions are present the definition of solution points brings along the fact the many ``non-dominated'' (equivalent) solutions exist. Such solutions compose the so-called Pareto front and there is great interest in trying to generate a set of points that approximate the Pareto front.

Hence, in this paper,  we are interested in developing new globally convergent derivative-free methods with posteriori articulation of preferences 
for nonlinearly constrained multiobjective minimization problems of the following form:
\begin{equation}\label{prob}
\begin{array}{ll}
 \min & F(x) = (f_1(x),\dots,f_q(x))^\top\\
 s.t. & g_i(x) \leq 0,\ \text{for all}\ i=1,\dots,m\\
 & l_i\leq x_i\leq u_i,\ \text{for all}\ i=1,\dots,n,\\
 & \ x_i\in \mathbb{R} \ \text{ for all } \ i\in I^c,
\\ [.3em] 
 &\ x_i\in \mathbb{Z} \ \text{ for all } \ i\in I^z,
\end{array}
\end{equation}
where $x\in  \mathbb{R}^n$, \, $l,u\in  \mathbb{R}^n$, \, and $I^c\cup I^z = \{1,{2,}\ldots ,n \}$, with $I^c \cap I^z = \emptyset$ and $I^c,I^z\neq\emptyset$. 
We assume $l_i < u_i$ for all $i = 1,\dots,n$, and
$l_i,u_i\in  \mathbb{Z}$ for all $i\in I^z$. 
$f_i:\Re^n\to \Re$, $i=1,\dots,q$, $g_i:\Re^n\to \Re$, $i=1,\dots,m$. 
We denote by $X$ and $Z$ the following sets
$$
X=\{x\in \Re^n: l_i\le x_i\le u_i,\ i=1,\dots,n\},\quad Z = \{x\in\Re^n:x_i\in\mathbb{Z}, i\in I^z\},
$$
and by $\cal F$ the feasible set of problem (\ref{prob}), namely,
$$
{\cal F} = \{x\in \Re^n: g_i(x) \leq 0,\ i=1,\dots,m\}\cap X\cap { Z}.
$$
We note that, by definition, $X$ is a compact set.\\

To solve problem \ref{prob} many methods have been proposed in the past years. These methods belong both to the class of deterministic and probabilistic methods.

Within the class of deterministic derivative-free methods for multiobjective mixed-integer problems, {\tt BiMADS} \cite{audet2008multiobjective} manages bi-objective problems and is able to generate an approximation of the Pareto front by solving a series of single-objective formulations of the bi-objective problem. The single objective problems are solved by using the MADS algorithm for single objective optimization. {\tt BiMADS} is able to manage also the presence of discrete or integer variables. The Directional Multi Search (DMS) algorithm proposed in  \cite{custodio2011direct} (and successively studied and enhanced in \cite{custodio2021worst,bras2020use,custodio2012recent}) was the first method in which the iterate is defined by a set o points rather than a single point. DMS can handle multiobjective problems with more than two objective functions but it cannot be used when some of the variables must take integer or discrete values. In \cite{liuzzi2016derivative} the DFMO algorithm has been proposed. DFMO is a linesearch-based derivative-free algorithm which evolves a set of non-dominated points that approximate Pareto stationary points in the limit. \cite{cocchi2018implicit} proposes {\tt imfil}, an adaptation to multiobjective derivative-free problems of the well-known implicit filtering algorithm for single objective. The {\tt imfil} method can be heuristically modified in such a way to produce a set of non-dominated points. 
More recently, in \cite{bigeon2021dmulti} a new MADS-type algorithm has been proposed for multiobjective problems, namely {\tt DMulti-MADS}. This method does not aggregate objectives and keeps a list of non dominated points which converges to a (local) Pareto set.

As concerns methods that try to approximate the set of global Pareto solutions, {\tt MODIR} \cite{campana2018multi} proposes and hybridization of the DIRECT-type algorithm {\tt DIRMIN} for derivative-free global optimization \cite{liuzzi2016exploiting} with the aforementioned DFMO \cite{liuzzi2016derivative}. In \cite{custodio2018multiglods} an algorithm is presented that alternates between initializing new searches, using a clever multistart strategy, and exploring promising subregions, by means of a directional direct search algorithm. Finally, in this context, the paper \cite{niebling2019branch} proposes a new branch-and-bound-based algorithm for smooth nonconvex multiobjective optimization problems with convex constraints. The algorithm computes an approximation of the set of Pareto global optimal solutions.

Probabilistic (and heuristic) approaches,  like e.g. genetic algorithms \cite{deb2002fast} and simulated annealing \cite{suppapitnarm2000simulated}.








\par\medskip


Inspired by the ideas in \cite{fasano2014linesearch} and \cite{custodio2011direct}, 
we describe a new exact-penalty-based linesearch approach (with sufficient decrease) for nonlinearly constrained MOO problems. 
In order to study the convergence properties of the proposed method, we  carry out a preliminary theoretical analysis of the problem itself. 
We describe new optimality conditions that take explicitly into account the bound constraints and 
that are obtained by only assuming Lipschitz continuity of the problem functions. We also prove that the original problem 
is equivalent to a bound constrained problem obtained by penalizing the nonlinear constraints with an exact merit function. 
In particular, we introduce a merit function that penalizes the general nonlinear inequality constraints in each term of $F(x)$ and we resort to the minimization of a penalty
function subject to the simple bound constraints. We would like to point out that some exact penalty methods for multiobjective optimization have already
been proposed in the literature (see e.g. \cite{fukuda2015external,huang2002nonlinear} and references therein). Anyway, to the best of our knowledge, 
this is the first time that a penalty approach with explicit handling 
of the bound constraints for Lipschitz continuous multiobjective problems is studied.
The proposed approach enables to handle also infeasible starting points. Furthermore, thanks to this exact penalty, we can adapt the derivative-free approach in \cite{fasano2014linesearch} to the MOO case. 
This new approach gives us three relevant advantages:\\
\begin{enumerate}
 \item[-]  by means of the sufficient decrease we can avoid the use of integer lattices 
 (see e.g. \cite{custodio2011direct, fasano2014linesearch});\\
 \item[-]  the extrapolation phase allows us to better exploit a descent direction and hence to prove, under some density assumptions on the search
directions, convergence to a set 
of Pareto stationary points (i.e. we prove that any accumulation point of the sequences generated by our method
is a Pareto stationary point);\\
 \item[-]  thanks to the exact penalty approach, the starting point can be infeasible.\\
\end{enumerate}

The fact that any accumulation point of the generated sequences is a Pareto stationary point, is an interesting theoretical result, since it 
is slightly stronger than the results reported in  \cite{custodio2011direct}. 
We would also like to note that, to the best of our knowledge, the idea of only penalizing the nonlinear constraints is new in the context of multiobjective derivative-free optimization. 

In the last part of the paper, we test the linesearch approach on both bound constrained and nonlinearly constrained problems. The aim of the tests on 
bound constrained problems is  understanding to what extent the theoretical properties of our method are helpful in practice.
On the other hand, the goal of the tests on nonlinearly constrained problems is to show the effectiveness of the exact penalty approach 
when embedded in a DFO method for bound constrained multiobjective problems. 
For this reason we report the results obtained by both our algorithm and the globally convergent version of DMS.

\subsection{Notations and preliminary material}\label{sec:notation}
We start the section by adapting the Lipschitz continuity to the mixed integer case. 
\begin{definition}[Lipschitz continuity w.r.t. the continuous variables]\label{def:Lipshitz}
Let $h: \mathbb{R}^n \to \mathbb{R}$. We say that $h$ is Lipschitz continuous w.r.t. the continuous variables when $L_h$ exists such that, for all $x,y\in\Re^n$ with $x_i=y_i$, for $i\in I^z$, it results
\[
 \|h(x)-h(y)\| \leq L_{h}\|x-y\|. 
\]
\end{definition}
Then, we state the main assumption on the functions defining problem (\ref{prob}).
\begin{ass}[Lipschitz continuity]\label{ass:Lipscitz}
We assume that all the functions defining problem (\ref{prob}), namely $f_i$, $i=1,\dots,q$, and $g_j$, $j=1,\dots,m$, are Lipschitz continuous w.r.t. the continuous variables (see definition \ref{def:Lipshitz})
\end{ass}
\begin{definition}[Generalized directional derivative and generalized gradient] \label{def:clarke_directional_derivative}
	Let $h: \mathbb{R}^n \to \mathbb{R}$ be a Lipschitz continuous function near $x \in \mathbb{R}^n$ with respect to its continuous variables $x_c$ (see definition \ref{def:Lipshitz}). The generalized directional	derivative of $h$ at $x$ in the direction $s \in \mathbb{R}^n$, with $s_i = 0$ for $i \in I^z$, is 
	\begin{equation}\label{eq:clarke_directional_derivative}
		h^{Cl}_{{{c}}}(x; s) = \limsup_{\footnotesize\begin{array}{l}y_c\to x_c, y_z = x_z, t\downarrow 0\end{array}} \frac{h(y+t s) -
			h(y)}{t}.
	\end{equation}
	To simplify the notation, the generalized gradient of $h$ at $x$ w.r.t the continuous variables can be redefined as
	\begin{eqnarray*}
		\partial_{{{c}}} h( x) & = & \{v\in \mathbb{R}^{n} :\ v_i=0,\ i\in I^z,\ \mbox{and}\  h_{{{c}}}^{Cl}( x; s)\ge s^Tv \ \text{ for all } s \in \mathbb{R}^n, \\
		&& \quad \text{ with } s_i = 0 \text{ for } i \in I^z \}.
	\end{eqnarray*}
\end{definition}
\smallskip
Given a vector $v\in\Re^n$, a subscript will be used to denote either one of its components ($v_i$) or the fact that it is an
element of an infinite sequence of vectors ($v_k$). In case of possible misunderstanding or ambiguities, the $i$th component of a
vector will be denoted by $(v)_i$. 

Given two vectors $u,v\in\Re^n$, we use the following convention for vector equalities and inequalities:
\begin{eqnarray*}
 u = v & \Leftrightarrow & u_i = v_i,\ \forall\ i=1,\dots,n,\\
 u < v     & \Leftrightarrow & u_i < v_i,\ \forall\ i=1,\dots,n,\\
 u \leq v  & \Leftrightarrow & u_i \leq v_i,\ \forall\ i=1,\dots,n,\ \mbox{and}\ u\neq v,\\
 u \leqq v & \Leftrightarrow &  u_i \leq v_i,\ \forall\ i=1,\dots,n.
\end{eqnarray*}
Note that $u \geq v$ if and only if $-u \leq -v$.

We denote by $v^j$ the generic $j$th element of a finite set of vectors. Furthermore, vectors $e_1,\dots,e_n$ represent the coordinate unit vectors.
Given two vectors
$a,b\in\Re^n$, we denote by $y=\max\{a,b\}$ the vector such that $y_i=\max\{a_i,b_i\}$, $i=1,\dots,n$. Furthermore, given a
vector $v$ we denote by $v^+=\max\{0,v\}$. The projection of a point $x$ onto the set $X$ will be denoted by $[x]_{[l,u]}$.
Finally,  we denote the unit sphere in the origin by $S(0,1) = \{d\in\Re^n: \|d\|=1\}$, and $Co(A)$ indicates
the convex hull of the set $A$. Given a point $x$ and a scalar $\rho > 0$, ${\cal B}(x,\rho) =
\{y\in\Re^n: \|x-y\|\leq \rho\}$. We denote 
$$\Gamma = \{\mu\in\Re^q:\mu\geq 0, \sum_{i=1}^q\mu_i = 1\}.$$ 
Finally, by $\one$ we denote the vector of all ones of dimension $q$.

When dealing with several objective functions at a time, the concept of Pareto dominance is usually considered in the comparison of two points.

\begin{definition}[Pareto dominance]
Given two points, $x,y\in\cal F$, we say that $x$ dominates $y$ if $F(x)\leq F(y)$.
\end{definition}
\par\smallskip\noindent

Anyway, when coming to optimality, it may not be possible to find a point which is optimal for
all the objectives simultaneously. This is the reason why the concept of Pareto dominance is also used
to characterize global and local optimality into a multiobjective framework. More specifically, by means of the following
two definitions, we are able to identify a set of nondominated points (the so called Pareto front or
frontier) which represents the (global or local) optimal solutions of a given multiobjective problem.

\begin{definition}[Global Pareto optimality]
A point $x^\star\in\cal F$ is a global Pareto optimizer of Problem (\ref{prob}) if it does not exist a point $y\in\cal F$ such
that $F(y)\leq F(x^\star)$.
\end{definition}


Then, we have to define local Pareto solutions. Hence, we have to formally state the definition of neighborhoods in the mixed integer contexts. To this end, with respect to the discrete variables, we have to preliminary introduce the concept of feasible primitive directions which will be used to define discrete neighborhoods.

\begin{definition}[set of feasible primitive directions]
	Given $x\in X\cap Z$,
	\vskip -0.3cm
	\[
	D^z(x) = \{d\in\mathbb{Z}^n:\ d_z\ \mbox{is primitive\footnote{the greatest common divisor of its component is 1}}, d_c = \mathbf{0}, x+ d\in X\cap Z\}.
	\] 
	\vskip -0.1cm
	It can be proved that $\bar D = \bigcup_{x\in X\cap Z} D^z(x)$ is finite
\end{definition}

\begin{definition}[Neighborhoods]	
	Given $\bar x\in X\cap Z$ and $\rho > 0$,
	\vskip -0.7cm
	\begin{eqnarray*}
		{\cal B}^z(\bar x) & = & \{x\in X\cap Z:\ x = \bar x+d,\ \forall\ d\in D^z(\bar x)\},\\
		{\cal B}^c(\bar x;\rho) & = & \{x\in X\phantom{\cap Z}\ :\ x_z=\bar x_z,\ \|x_c-\bar x_c\|\leq \rho\}
	\end{eqnarray*}
\end{definition}

Finally, we can give the definition of local Pareto solutions for problem (\ref{prob}).
\begin{definition}[local Pareto optimum]\label{def:locpareto}
	$x^*\in X\cap Z$ is a local Pareto optimum of problem (\ref{prob}) if
	\begin{enumerate}
		\item there is no $y\in {\cal B}^z(x^*)\phantom{,\rho}\cap{\cal F}$ s.t. $F(y) \leq F(x^*)$ 
		\item there is no $y\in {\cal B}^c(x^*,\rho)\cap{\cal F}$ s.t. $F(y) \leq F(x^*)$ for some $\rho > 0$.
	\end{enumerate}
\end{definition}	

As usual in the continuous optimization setting, we have now to introduce the definition of ``stationary'' points. For this reason, and to explicitly handle the simple bound constraints, we first need to formally define the cone of feasible directions for the bound constraints and w.r.t. the continuous variables.

\begin{definition}[Set of  feasible continuous directions]\label{def:cone_feasible_continuous_directions} Given a point $x\in X\cap\cal{Z}$, the set
	\begin{eqnarray*}
		D^c(x) & = &\{s\in \mathbb{R}^n: s_i=0\quad \hbox{for all}\quad i\in I^z,\\
		&   &\hspace*{1.5cm} s_i\ge 0 \quad \hbox{for all}\quad  i\in I^c\quad\mbox{and}\quad x_i = l_i,\\
		& &\hspace*{1.5cm}  s_i\le 0\quad \hbox{for all}\quad  i\in I^c\quad\mbox{and}\quad  x_i = u_i, \\ 
		&   & \hspace*{1.5cm} s_i\in \mathbb{R}\quad \hbox{for all}\quad  i\in I^c\quad\mbox{and}\quad l_i<x_i < u_i,\}
	\end{eqnarray*}
	is the cone of feasible continuous directions at $x$ with respect to $X\cap\cal{Z}$.
\end{definition}

Finally, stationary points of problem (\ref{prob}) can be defined.

\begin{definition}[Pareto-Clarke Stationarity w.r.t. continuous variables]\label{Pareto_Clarke_stationarity}
A point $x^*\in\cal F$ is a Pareto-Clarke stationary point w.r.t. the continuous variables when there exist multipliers
	$\sigma_1^\star,\dots,\sigma_q^\star,\lambda_1^\star,\dots,\lambda_m^\star\in\Re$, such that 
	\[
	\begin{array}{l}
		\qquad \sigma_i^\star\geq 0,\quad \sigma^\star \neq 0,\qquad i=1,\dots,q,\\
		\qquad \lambda_j^\star\ge 0,\quad \lambda_j^\star g_j(x^\star)=0,\qquad j=1,\dots,m,
	\end{array}
	\]
	and a vector 
	\[
	\begin{array}{l}
		\qquad \bar \xi\in \displaystyle\sum_{i=1}^q\sigma_i^\star \partial_c f_i(x^\star) + \displaystyle\sum^m_{j=1}\lambda_j^\star\partial_c g_j(x^\star),
	\end{array}
	\]
	such that
	\[
	\bar \xi^\top d\geq 0, \quad \mbox{for every $d\in D^c(x^\star)$}.
	\]
\end{definition}

Then, we say that a point $x^*\in\cal F$ is (Pareto) stationary for problem \ref{prob} when:
\begin{itemize}
    \item[(i)] it is Pareto-Clarke stationary w.r.t. the continuous variables and
    \item[(ii)] it is a local Pareto optimum w.r.t. the discrete variables.
\end{itemize}

\begin{proposition}[Pareto-Clarke KKT Necessary Optimality Conditions]\label{KKT_cond}
	
	Let $x^\star\in\cal F$ be a local Pareto minimum of the problem \eqref{prob} and assume that a direction $d\in D^c(x^\star)$ exists such that
	for all $ j\in\{1,\dots,m: g_j(x^\star) = 0\}$ :
	
	\begin{eqnarray}
		(\xi^{g_j})^\top d < 0, \quad \mbox{for all}\quad \xi^{g_j}\in\partial_c g_j(x^\star).\label{mfcq2}
	\end{eqnarray}
	Then, $x^*$ is a Pareto stationary point.
\end{proposition}

{\bf Proof}. The proof easily follows recalling definition \ref{def:locpareto} and the proof of \cite[Proposition 3.4]{liuzzi2016derivative}.
$\hfill\Box$
\par\medskip

\section{The bound constrained case}

\begin{equation}\label{probbox}
	\begin{array}{ll}
		\min & F(x) = (f_1(x),\dots,f_q(x))^\top\\
		s.t. & x\in X\cap Z,
	\end{array}
\end{equation}


\begin{proposition}\label{st_pt2}
	A point $\bar x\in X\cap Z$ is a  Pareto-Clarke stationary point of (\ref{probbox}) w.r.t. the continuous variables,
	if and only if  for all $d\in D^c(\bar x)$, an index
	$j_d\in\{1,\dots,q\}$ exists such that:
	\begin{equation*}
		(f_{j_d})_c^{Cl}(\bar x;d) \geq 0.
	\end{equation*}
\end{proposition}
\par\smallskip\noindent
{\bf Proof}. The proof is an almost trivial repetition of \cite[Proposition 3.7]{liuzzi2016derivative}.$\hfill\Box$
\par\medskip\noindent

We remark that a point $\bar x\in X\cap Z$ is (Pareto) stationary for problem (\ref{probbox}) if and only if:
\begin{itemize}
    \item[(i)] it is Pareto-Clarke stationary  w.r.t. the continuous variables and
    \item[(ii)] it is a local Pareto optimum w.r.t. the discrete variables. 
\end{itemize}

Given a point $x \in X\cap Z$, the Clarke-Jahn generalized
directional derivative of a function $f$ along the direction $d \in D^c(x)$ is given by (see \cite[Section 3.5]{Jahn:96}):
\begin{equation}\label{clarke-jahndd}
f^\circ(x; d) = \limsup_{\footnotesize\begin{array}{l}y_c\to x_c, y_z=x_z, y\in X\cap Z\\ t\downarrow 0, y+t d\in X\cap Z\end{array}} \frac{f(y+t d) -
f(y)}{t}.
\end{equation}
Recall that, according to definitions \eqref{eq:clarke_directional_derivative} and \eqref{clarke-jahndd}, we have, for any $d\in\Re^n$,
\begin{equation}\label{ineqCJ}
 f^{Cl}(x;d)  \geq f^\circ(x;d).
\end{equation}

\begin{definition}\label{def_pareto-clarke-jahn-staz} A point $x^*\in X\cap Z$ is a Pareto-Clarke-Jahn stationary point w.r.t. the continuous variables if, for all $d\in D^c(x^*)$, an index $j_d\in\{1,\dots,q\}$ exists such that
\[
  (f_{j_d})_c^\circ(\bar x; d) \geq 0.
\]
\end{definition}

\begin{lemma}\label{rel_jahn-clarke}
With respect to the continuous variables, every Pareto-Clarke-Jahn stationary point $x^*\in X\cap Z$ is also a Pareto-Clarke stationary point.
\end{lemma}
{\bf Proof}. 
Recalling the definitions (\ref{eq:clarke_directional_derivative}) and (\ref{clarke-jahndd}) of Clarke and Clarke-Jahn generalized directional derivatives, respectively, and relation (\ref{ineqCJ}), we have that
\[
   (f_{j_d})_c^{Cl}(\bar x;d) \geq (f_{j_d})_c^\circ(\bar x;d).
\]
Hence, the result easily follows.
$\hfill\Box$

\section{The penalty approach}

Given problem (\ref{prob}), we introduce the following nonsmooth penalty functions
\[
Z_i(x;\eps) = f_i(x) + \frac{1}{\eps}\sum_{j=1}^m\max\{0,g_j(x)\}
\]
and consider the bound constrained problem
\begin{equation}\label{prob_Z}
\begin{array}{ll}
 \min & Z(x;\eps) = (Z_1(x;\eps),\dots,Z_q(x;\eps))^\top\\
 s.t. & x\in X\cap Z
\end{array}
\end{equation}

\begin{ass}[EMFCQ for mixed-integer problems]\label{assmfcq}
	Given Problem \eqref{prob}, for any \\ $x\in (X\cap{Z})\setminus\stackrel{\circ}{\cal F}$, one of the following conditions holds:
	\begin{itemize}
		\item[(i)] a direction $s\in D^c(x)$ exists such that  \[
		(\xi^{g_i})^\top s < 0,
		\]
		for all $\xi^{g_i}\in\partial_c g_i(x)$ with $i\in\{h \in \{1,\dots,m\}: \ g_h(x)\geq 0\}$; 
		\item[(ii)] a direction $\bar d \in D^z(x)$ exists such that
		$$\sum_{i=1}^m \max\{0, g_i(x+\bar d)\} < \sum_{i=1}^m \max\{0, g_i(x)\}.$$ 
	\end{itemize}
\end{ass}

We first prove the equivalence between local minimum points and global minimum points of the original problem (\ref{prob}) and the penalized problem (\ref{prob_Z}). 


\begin{proposition}\label{th:no_stationary_points_out_feasible_region}	
	Let Assumption \ref{assmfcq} hold. A threshold value $\varepsilon^\star > 0$ exists such that the function $P(x;\varepsilon)$ has no Pareto-Clarke stationary points in $(X\cap\mathcal{Z})\setminus\cal F$ for any $\varepsilon \in (0,\varepsilon^\star]$. 
\end{proposition}
{\bf Proof}. \quad We proceed by contradiction and assume that for any integer $k$, an $\varepsilon_k \leq 1/k$ and a stationary point for Problem \eqref{prob_Z} $x_k\in (X\cap\mathcal{Z})\setminus\cal F$ exist. Then, let us consider a limit point $\bar x\in (X\cap\mathcal{Z})\setminus\stackrel{\circ}{\cal F}$ of the sequence $\{x_k\}$ and, without loss of generality, let us call the corresponding subsequence as $\{x_k\}$ as well. Then, since $x_k\to\bar x$, the discrete variables remain fixed, i.e. $(x_k)_i = \bar x_i$, for all $i\in I^z$ and $k$ sufficiently large. Now, the proof continues by separately assuming that point $(i)$ or $(ii)$ of Assumption \ref{assmfcq} holds. 
\par
First we assume that point $(i)$ of Assumption \ref{assmfcq} holds at $\bar x$.
Therefore, a direction $\bar s\in D^c(\bar x)$ exists such that 
\[
\left(\xi^{g_i}\right)^\top \bar s <0 \ \mbox{ for all }\ \xi^{g_i}\in\partial_c g_i(\bar x) \mbox{ with } i\in\{ h \in \{1,\dots,m\}: g_h(\bar x) \geq 0\}.
\]

In particular, it follows that
\begin{equation}\label{maxxi}
\max_{\tiny\begin{array}{c}\xi^{g_i}\in\partial_c g_i(\bar x)\\ i\in I(\bar x)\end{array}} \left(\xi^{g_i}\right)^\top \bar s = -\eta < 0,
\end{equation}
where $I(\bar x) = \left\{i\in\{1,\dots,m\}: g_i(\bar x) = \bar \phi(\bar x) \right\}$, $\bar \phi(x) = \max\left\{0,g_1(x),\dots,g_m(x)\right\}$ and $\eta$ is a positive scalar. Note that $\bar x\not\in\cal F$ implies $\bar \phi(\bar x) > 0$. 
\par
By \cite[Proposition 2.3]{lin2009decomposition}, it follows that $\bar s\in D^c(x_k)$. Moreover, since $x_k$ is a stationary point, we have 
that an index
	$j\in\{1,\dots,q\}$ exists such that:
	\begin{equation}\label{fzero_cl}
	\max_{\xi\in\partial_c Z_j(x_k;\varepsilon)}\xi^\top \bar s = 	(P_{j})_c^{Cl}(x_k;\eps,\bar s) \geq 0.
	\end{equation}
 	and there is no $y\in {\cal B}^z(x_k)$ s.t. $F(y) \leq F(x_k)$.

By \cite{clarke1990optimization}, we know that 
\begin{equation}
\label{eq:partial_1}
\partial_c Z_j(x_k;\varepsilon) \subseteq \partial_c f_j(x_k) + \frac{1}{\varepsilon}\partial_c (\max\left\{0,g_1(x_k),\dots,g_m(x_k)\right\})
\end{equation}
and 
\begin{equation}
\label{eq:partial_2}
\partial (\max\left\{0,g_1(x_k),\dots,g_m(x_k)\right\}) \subseteq Co\left(\{\partial_c g_i(x_k): i\in I(x_k)\}\right),
\end{equation}
where 
$Co(A)$ denotes the convex hull of a set $A$ (see \cite[Theorem 3.3]{Rockafellar.1970}).
\par
Therefore, by \eqref{fzero_cl}--\eqref{eq:partial_2}, $\xi^{f}_k \in \partial_c f_j(x_k)$, \, $\xi^{g_i}_k \in \partial_c g_i(x_k)$ and $\beta^i_k$ \mbox{ with } $i\in I(x_k)$ exist such that
\begin{equation}\label{derposxi}
\left(\xi^{f}_k + \frac{1}{\varepsilon_k}\sum_{i\in I(x_k)}\beta^i_k\xi^{g_i}_k\right)^\top \bar s \geq 0,
\end{equation}

$$
\sum_{i\in I(x_k)}\beta^i_k = 1 \ \text{ and } \ \beta^i_k \geq 0.
$$
Since $m$ is a finite number, there exists a subsequence of $\{x_k\}$ such that $I(x_k) = \bar I$. 

Then, recalling that $(x_k)_i=\bar x_i$, for all $i\in I^z$ and for $k$ sufficiently large, and since a locally Lipschitz continuous function has a generalized gradient which is locally bounded, it results that the sequences $\{\xi^{f}_k\}$ and $\{\xi^{g_i}_k\}$, with $i\in \bar I$, are bounded. Hence, we get that
\begin{subequations}\label{limconds}
	\begin{eqnarray}
	&& \label{limconds1}\xi^{f}_k\to \bar\xi^{f},\\
	&& \xi^{g_i}_k\to \bar\xi^{g_i} \ \mbox{for all}\ i\in \bar I,\\
	&& \beta^i_k\to\bar\beta^i \ \ \mbox{for all}\ i\in \bar I.
	\end{eqnarray}
\end{subequations}
Now the upper semicontinuity of $\partial_c f_j$ and $\partial_c g_i$, with $i\in\bar I$, at $\bar x$ (see Proposition 2.1.5
in \cite{Clarke83}) implies that $\bar\xi^{f}\in \partial_c f_j(\bar x)$ and $\bar\xi^{g_i}\in \partial_c g_i(\bar x)$ for all
$i\in\bar I$.

The continuity of the problem functions guarantees that for $k$ sufficiently large
\[
\{i:g_i(\bar x)-\phi(\bar x) < 0\} \subseteq \{i: g_i(x_k) -\phi(x_k) < 0\}, 
\]
and, in turn, this implies that for $k$ sufficiently large
\[
\{i: g_i(x_k) -\phi(x_k) = 0\} = I(x_k) \subseteq I(\bar x) = \{i:g_i(\bar x)-\phi(\bar x) = 0\}. 
\]
Since $I(x_k) \subseteq I(\bar x)$, we have that
\begin{equation}\label{inclusione}
\bar I \subseteq I(\bar x). 
\end{equation}
Finally, for $k$ sufficiently large, \eqref{maxxi}, \eqref{limconds}, and \eqref{inclusione} imply 
\begin{equation}\label{derposxi2}
\left(\xi_k^{g_i}\right)^\top\bar s \leq -\frac{\eta}{2} \ \text{ for all } i\in\bar I,
\end{equation}
and (\ref{derposxi}), multiplied by $\varepsilon_k$, implies
\begin{equation} \label{eq:eq_3}
\left(\varepsilon_k\xi^{f}_k +\sum_{i\in\bar I}\beta^i_k\xi^{g_i}_k\right)^\top \bar s \geq 0.
\end{equation}
Equations \eqref{derposxi2} and \eqref{eq:eq_3} yield
\[
0 \leq \left(\varepsilon_k\xi^{f}_k +\sum_{i\in\bar I}\beta^i_k\xi^{g_i}_k\right)^\top \bar s \leq 
\left(\varepsilon_k\xi^{f}_k \right)^\top \bar s -\frac{\eta}{2}.
\]
which, by using \eqref{limconds}, gives rise to a contradiction when $\varepsilon_k \to 0$.

Now we assume that point $(ii)$ of Assumption \ref{assmfcq} holds at $\bar x$. Let  $\bar d \in D^z(\bar x)$ be the direction such that
\begin{equation}
	\label{eq:eq_dir}
	\sum_{i=1}^m \max\{0, g_i(\bar x+\bar d)\} < \sum_{i=1}^m \max\{0, g_i(\bar x)\}.
\end{equation} 
recalling that $(x_k)_i=\bar x_i$, for all $i\in I^z$ and for $k$ sufficiently large, we have that for $k$ sufficiently large $D^z(\bar x) = D^z(x_k)$, so that $\bar d\in D^z(x_k)$. 
By definition of stationary point and discrete neighborhood, we have that either
$$
Z(x_k; \varepsilon_k)\leqq Z(x_k + \bar d; \varepsilon_k),$$
or an index $j$ exists such that
$$Z_j(x_k; \varepsilon_k) < Z_j(x_k + \bar d; \varepsilon_k),$$
{where} $x_k + \bar d \in {\cal B}^z(x_k)$.

Hence, an index $h$ exists such that
$$f_h(x_k) + \frac{1}{\varepsilon_k} \sum_{i=1}^m \max\{0, g_i(x_k)\} \le f_h(x_k + \bar d) + \frac{1}{\varepsilon_k} \sum_{i=1}^m \max\{0, g_i(x_k+\bar d)\}.$$
Multiplying by $\varepsilon$ and considering that $\varepsilon_k \to 0$, if we take the limit for $k \to \infty$ we have that
$$\sum_{i=1}^m \max\{0, g_i(\bar x)\} \le \sum_{i=1}^m \max\{0, g_i(\bar x+\bar d)\}.$$
The latter equation is in contradiction with \eqref{eq:eq_dir}.
$\hfill\Box$ \\
\par
Now, we can prove that there exists
a threshold value $\bar \varepsilon$ for the penalty parameter such that, for any  $\varepsilon \in(0, \bar \varepsilon)$, any local minimum of the penalized problem is also a local minimum of the original problem.

\par\medskip\noindent

\begin{proposition}\label{equivstat}
For any $\eps >0$, every Pareto-Clarke stationary point $\bar x$ of Problem (\ref{prob_Z}), such that $\bar x\in{\cal F}$, is also a
Pareto-Clarke stationary point of Problem (\ref{prob}).
\end{proposition}
\par\smallskip\noindent
{\bf Proof}. Since $\bar x$ is, by assumption, a Pareto-Clarke stationary point of Problem (\ref{prob_Z}), then by Definition \ref{Pareto_Clarke_stationarity} we know that a
vector of non-negative multipliers $\mu\in\Gamma$, i.e., $\mu\geq 0$, $\sum_{i=1}^q\mu_i$, not all zero,  and a vector  
$\xi^\star\in\sum_{i=1}^q\mu_i\partial_c Z_i(\bar x,\eps)$ exist such that, for all $d\in D^c(\bar x)$,
$$(\xi^\star)^\top d \geq 0.$$
Now, we recall that
$$
 \partial_c Z_i(x,\eps) \subseteq \partial_c f_i(x) + \frac{1}{\eps}\sum_{j\in I(x)} \partial_c g_j(x),
$$
where  $I(x) = \{i:\ g_i(x)\geq 0\}$.
Hence, we have that $\xi^* \in \sum_{i=1}^q\mu_i \left ( \partial_c f_i(\bar x) + \frac{1}{\eps}\sum_{j\in I(\bar x)} \partial_c g_j(\bar x)\right )$, $\mu\in\Gamma$.
Then, denoting $\lambda_j = \left( \sum_{i=1}^q\mu_i \right )/\eps$,
\[
  \displaystyle\max \left \{ \xi^\top d : \ \xi\in \displaystyle\sum_{i=1}^q\mu_i\partial_c f_i(\bar x) + \displaystyle\sum_{j\in I(\bar x)} \lambda_j\partial_c g_j(\bar x) \right \}\geq 0
\]
for all $d\in D^c(\bar x)$ with $\lambda_j\geq 0$, $j\in I(\bar x)$. The above condition shows that a 
$\bar\xi\in \displaystyle\sum_{i=1}^q\mu_i\partial_c f_i(\bar x) + \displaystyle\sum_{j\in I(\bar x)} \lambda_j\partial_c g_j(\bar x)$ exists such that 
$\bar\xi^\top d \geq 0,\ \mbox{for all}\ d\in D^c(\bar x)$. Hence,
by recalling that $\bar x\in\cal F$ and $I(\bar x) = I_0(\bar x)$ when $\bar x\in\cal F$, and setting $\la_j=0$ when $j\not\in I_0(\bar x)$, 
it is proved that $\bar x$ is a Pareto-Clarke  stationary point of Problem \eqref{prob}.$\hfill\Box$
\par\medskip
Now, we introduce an intermediate result which basically states that, for $\eps$ sufficiently small, every Pareto-Clarke-Jahn stationary point of Problem \eqref{prob_Z}
is a Pareto-Clarke  stationary point of Problem (\ref{prob}).
\par\smallskip\noindent
\begin{proposition}\label{equivalenza}
Let Assumption \ref{assmfcq} hold. Given problem \eqref{prob} and considering problem \eqref{prob_Z}, a threshold value
$\eps^\star >0$ exists such that, for every $\eps\in (0,\eps^\star]$, every Pareto-Clarke-Jahn-stationary point of Problem \eqref{prob_Z}
is a Pareto-Clarke  stationary point of Problem (\ref{prob}).
\end{proposition}
\par\smallskip\noindent
{\bf Proof}. Let $\bar x\in X\cap Z$ be Pareto-Clarke-Jahn stationary for Problem \eqref{prob_Z}. By \eqref{ineqCJ},
we also have that $\bar x$ is Pareto-Clarke stationary for Problem \eqref{prob_Z}.
Now, the proof follows by considering Propositions \ref{th:no_stationary_points_out_feasible_region} and \ref{equivstat}. $\hfill\Box$
\par\medskip\noindent

\par\medskip

We conclude by reporting from \cite{liuzzi2016derivative} a result concerning correspondence between global Pareto minimizers of problems (\ref{prob_Z}) and (\ref{prob}).

\begin{proposition}[{\cite[Proposition 3.13]{liuzzi2016derivative}}]\label{equivglobal}
Let Assumption \ref{assmfcq} hold. Then, given Problem \eqref{prob} and considering Problem \eqref{prob_Z}, a threshold value $\eps^\star >0$ exists such that, for every
$\eps\in (0,\eps^\star]$, any global Pareto minimizer of Problem \eqref{prob_Z} is a global Pareto minimizer of Problem \eqref{prob},
and conversely.
\end{proposition}

\section{The algorithm}
In the definition of Algorithm DFMOINT, and in particular in Phase 2.a (Explore discrete variables), we shall make use of the following notation. Given a list 
\[
  L = \{(x_j,\alpha_j^c,\alpha_j^{(d)},\xi_j):\ j=1,2,\dots,|L|\}
\]
we denote 
\[
  L[x] = \{x_j: (x_j,\alpha_j^c,\alpha_j^{(d)},\xi_j)\in L\}.
\]
Note that, $L$ is a set of tuples whereas $L[x]\subset\Re^n$ is a set of points in $\Re^n$.

\begin{algorithm}
	\caption{
		DFMOINT
	}\label{alg_mixed_3}
	\begin{algorithmic}[1]
		\par\medskip
		\item \ Let $x_0\in X\cap\mathcal{Z}$, \ $\xi_0>0$, \  $\theta\in (0,1)$,\ $\eps > 0$;
		\smallskip
		\item \ let $\{s_k\}$ be a sequence such that $s_k\in S(0,1)$ for all $k$; let $D=D_0$ be a set of initial primitive discrete directions;
		\smallskip
		\item \ let $\tilde\alpha_0^{c} = 1$ be the initial stepsize along $s_k$; let $\tilde\alpha_0^{(d)} = 1$ be the initial stepsizes along  $d \in D$;
		\smallskip
		\item \ let $L_0= \{(x_i,\alpha_i^c,\alpha_i^{(d)},\xi_i):x_i\in X\cap Z,\alpha_i^c = \tilde\alpha_0^c,\alpha_i^{(d)} = \tilde\alpha_0^{(d)}, i=1,\dots,|L_0|\}$

		\par\medskip
		
		\item {\bf For} $k=0,1,\dots$
		\par\medskip
		\item \qquad Set $\tilde L_k = L_k$
		
		\Algphase{\qquad\quad\fbox{PHASE 1 - Explore continuous variables }}
		\smallskip		\item \qquad {\bf For} $i=1,2,\dots,|L_k|$
		\smallskip
		\item \qquad\quad Let $(x_i,\alpha_i^c,\alpha_i^{(d)},\xi_i)$ be the $i$-th tuple of $L_k$ 
		\smallskip

		\item \qquad\quad Compute $\tilde L_k$ by the {\em Projected Expansion}$((x_i,\alpha_i^c,\alpha_i^{(d)},\xi_i),s_k,\tilde L_k)$.
		\item \qquad {\bf End For}
 
		\item \qquad Set $L^c_k=\tilde L_k$ and  $L_k^d = L_k^c$
		\Algphase{\qquad\quad\fbox{PHASE 2.A - Explore discrete variables}}
		\smallskip
        \item \qquad {\bf For} $j=1,2,\dots,|L_k^c|$ 
		\par\medskip
		\item\qquad\quad Let $(x_j,\alpha_j^c,\alpha_j^{(d)},\xi_j)$ be the $j$-th tuple of $L_k^c$, set $D = D_k$ and $L^+ = L_k^d$
		\item \qquad \quad {\bf While} $D\neq \emptyset$ and $L^+[x]=L_k^d[x]$ {\bf do}
		\medskip
		
		
		\item \qquad\quad\quad Choose $d \in D$, set $D=D\setminus\{d\}$
		\item \qquad\quad\quad Compute $L_k^d$ by the {\em Discrete Search}$((x_j,\alpha_j^c,\alpha_j^{(d)},\xi_j),d,L_k^d)$.
		
		\medskip
		%

		\item \qquad\quad {\bf End While}
		\medskip
		\item \qquad\quad {\bf If} $L^+[x] = L_k^d[x]$ and $\alpha_j^{(d)} = 1$ for all $d\in D_k$ {\bf then}
		\item \label{step:red_xi}\qquad\quad\quad {\bf If} $(x_j,\alpha_j^c,\alpha_j^{(d)},\xi_j)\in L_k^d$ {\bf then }Set $L_k^d =  (L_k^d\setminus\{(x_j,\alpha_j^c,\alpha_j^{(d)},\xi_j)\})\cup\{(x_j,\alpha_j^c,\alpha_j^{(d)},\theta\xi_j)\}$
        \item \qquad\quad{\bf End If}
        \item \qquad{\bf End For}
		\par\smallskip

		\Algphase{\qquad\quad\fbox{PHASE 2.B - Update the set of primitive search directions}}
		\smallskip
		\item \qquad {\bf If} $L_k^c = L_k^d$  and  $\xi_i$ for all $(x_i,\alpha_i^c,\alpha_i^{(d)},\xi_i)\in L_k^d$ have been reduced at step \ref{step:red_xi} \textbf{then}
		
		\item \qquad\quad {\bf If} $ D_{k} =  \bar D^z$ {\bf then}
		
		\item \qquad\qquad\quad      Set $D_{k+1} = D_k$ and  $D=D_{k+1}$,

		\item \qquad\quad {\bf else}
		
		\item \qquad\qquad\quad Generate
		$D_{k+1}$  such that $D_{k+1} \subseteq \bar D^z$ and $ D_{k+1} \supset D_k$, set
		$D=D_{k+1}$.
		
		\item \qquad\qquad\quad Set $\alpha_{i}^{(d)} = 1$  for all $d\in D_{k+1}\setminus D_k$ and $(x_i,\alpha_i^c,\alpha_i^{(d)},\xi_i)\in L_k^d$.		
		
		\item \qquad\quad {\bf End If}
		
		\item \qquad {\bf else}
		\item \qquad\quad  Set $D_{k+1} = D_k$ and $D=D_{k+1}$.
		
		\item \qquad {\bf End If}
		\medskip
\smallskip
		\item \qquad Set $L_{k+1} = L_k^d$.
		\medskip
		\item {\bf End For}
		
		\par\bigskip\noindent
		
	\end{algorithmic}
\end{algorithm}

\begin{algorithm}
\caption{{{\tt Add\&Filter} ($L$, $(x,\alpha_x,\alpha^{(d)}_x,\xi_x)$)}}
\begin{algorithmic}
\par\medskip
 \item Set $\tilde L=\{(x,\alpha_x,\alpha^{(d)}_x,\xi_x)\}$
 \par\smallskip
 \item {\bf For each} $(y,\alpha_y,\alpha^{(d)}_y,\xi_y)\in L$\par
 \item \qquad{\bf If} $Z(x;\eps)\not\leq Z(y;\eps)$ {\bf then} $\tilde L = \tilde L\cup\{(y,\alpha_y,\alpha^{(d)}_y,\xi_y)\}$\par
 \item{\bf End For}\par
 \item Return $\tilde L$. 
\end{algorithmic}
\end{algorithm}

\par\medskip\noindent

\begin{algorithm}
	\caption{
		{\tt Projected Expansion} ($y,\hat\alpha,\alpha^{(d)},\xi,p,\tilde L$)
	}\label{alg:projexp}
	\begin{algorithmic}[1]
\par\smallskip
 \item $\delta,\theta \in (0,1)$, $\gamma > 0$.\par
 \item Set $\alpha = \hat\alpha$.\par
 \item \label{LS:step3a}{\bf If} $\exists\ x_j\in \tilde L$: $\Big(Z([y+\alpha p]_{[l,u]};\eps) > Z(x_j;\eps) - \gamma \displaystyle \alpha^2\one\Big)$ \  {\bf then}
 \item \qquad {\bf If} $(y,\hat\alpha,\alpha^{(d)},\xi)\in \tilde L$ {\bf then} $\tilde L = (\tilde L\setminus\{(y,\hat\alpha,\alpha^{(d)},\xi)\})\cup\{(y,\theta\hat\alpha,\alpha^{(d)},\xi)\}$
 \item \qquad Return $\tilde L$\par
 
 \item {\bf End If}
 
 \item\label{LS:step2} Let $\beta=\alpha/\delta$.\par
 \item \label{LS:step3}{\bf If} $\Big(Z([y+\beta p]_{[l,u]};\eps)\not < Z([y+\alpha p]_{[l,u]};\eps) - \gamma\left(\displaystyle \beta^2 -\alpha^2\right )\one\Big)$ {\bf then}\par
 
 \item \qquad \label{LS:step9} Set $\tilde L = \mbox{\texttt{Add\&Filter}}( \tilde L,([y+\alpha p]_{[l,u]},\alpha,\alpha^{(d)},\xi))$.
 
 \item {\bf Endif} 
 \item \label{LS:step4}{\bf If} $\Big(Z([y+\beta p]_{[l,u]};\eps)\not > Z(x_j;\eps) - \gamma \displaystyle \beta ^2\one\Big)$ \ $\forall\ x_j\in \tilde L$ {\bf then}\par 

 \item \qquad \label{LS:step12}Set $\alpha=\beta$ and go to Step \ref{LS:step2}.\par
  
 \item {\bf Endif}\par 
 
  \item Return $\tilde L$.

\end{algorithmic}
\end{algorithm}

\par\medskip\noindent

\begin{algorithm}
	\caption{
		{\tt Discrete Search} ($y,\tilde\alpha,\alpha^{(d)},\xi,p,\tilde L$)
	}\label{alg:DS}
	\begin{algorithmic}[1]
\par\smallskip
 \item $\delta \in (0,1)$, $\gamma > 0$.\par
 \item Compute the largest stepsize $\bar\alpha$ s.t. $ y+\bar\alpha p \in X\cap Z$
 \item Set $\alpha = \min\{\bar\alpha,\alpha^{(d)}_p\}$.\par
\item \label{DS:step4a} {\bf If} $\alpha = 0$ {\bf or} $\Big(Z(y+\alpha p;\eps) > Z(y;\eps) - \xi\one\Big)$ \  {\bf then}


 \item \qquad Let $\hat\alpha^{(d)} = \alpha^{(d)}$ and set $\hat\alpha^{(d)}_p =  \max\{1,\lfloor\alpha^{(d)}_p/2\rfloor$
  \item \qquad {\bf If} $(y,\tilde\alpha,\alpha^{(d)},\xi)\in\tilde L$ {\bf then} $\tilde L = (\tilde L\setminus\{(y,\tilde\alpha,\alpha^{(d)},\xi)\})\cup\{(y,\tilde\alpha,\hat\alpha^{(d)},\xi)\}$
  \item \qquad Return $\tilde L$\par
 
 \item {\bf End If}
 
 \item\label{DS:step2} Let $\beta=\min\{\bar\alpha,2\alpha^{(d)}_p\}$.\par
 \item \label{DS:step3}{\bf If} $\Big(Z(y+\beta p;\eps)\not < Z(y+\alpha p;\eps) - \xi\one\Big)$ {\bf then}\par
 
  \item \qquad Let $\hat\alpha^{(d)} = \alpha^{(d)}$ and set $\hat\alpha^{(d)}_p =  \alpha$
 
 \item \qquad \label{DS:step12} Set $\tilde L = \mbox{\texttt{Add\&Filter}}( \tilde L,(y+\alpha p,\tilde\alpha,\hat\alpha^{(d)},\xi))$.
 
 \item {\bf Endif} 
 \item \label{DS:step4}{\bf If} $\beta < \bar\alpha$ {\bf and} $\Big(Z(y+\beta p;\eps)\not > Z(x_j;\eps) - \xi\one\Big)$ \ $\forall\ x_j\in \tilde L$ {\bf then}\par 

 \item \qquad \label{DS:step15}Set $\alpha=\beta$ and go to Step \ref{DS:step2}.\par
  
 \item {\bf Endif}\par 
 
  \item {\bf If} $\beta = \bar\alpha$ {\bf then} 
  \item \qquad Let $\hat\alpha^{(d)} = \alpha^{(d)}$ and set $\hat\alpha^{(d)}_p =  \beta$
  \item \qquad \label{DS:step19} Set $\tilde L = \mbox{\texttt{Add\&Filter}}( \tilde L,(y+\beta p,\tilde\alpha,\hat\alpha^{(d)},\xi))$
  
  \item {\bf Endif}
  
  \item Return $\tilde L$.

\end{algorithmic}
\end{algorithm}

\par\medskip\noindent

\subsection{Convergence analysis}

\begin{definition}
Let $\{L_k\}$ with $L_k=\{(x_j,\alpha_j,\alpha_j^{(d)},\xi_j), \ j=1,\dots, |L_k|\}$ be the
sequence of sets of nondominated points produced by DFMOINT. We define a
\emph{linked sequence} as a sequence $\{(x_{j_k},\alpha_{j_k},\alpha_{j_k}^{(d)},\xi_{j_k})\}$ such that for any $k=1,2,\dots$, the pair
$(x_{j_k},\alpha_{j_k},\alpha_{j_k}^{(d)},\xi_{j_k})\in L_k$
is generated at iteration $k-1$ of DFMOINT by the pair $(x_{j_{k-1}},\alpha_{j_{k-1}},\alpha_{j_{k-1}}^{(d)},\xi_{j_{k-1}})\in L_{k-1}$.
\end{definition}

Let us, in particular, note the following cases.
\begin{enumerate}
    \item (Success w.r.t. continuous variables) we have
    \[
    \begin{split}
        & x_{j_k} = \left[x_{j_{k-1}} + \frac{\alpha_{j_{k-1}}}{\delta^{s_{k-1}}}d_{k-1}\right]_{[l,u]}, s_{k-1} \geq 0;\\
        & \alpha_{j_k} = \frac{\alpha_{j_{k-1}}}{\delta^{s_{k-1}}} \\
        & \alpha_{j_k}^{(d)} = \alpha_{j_{k-1}}^{(d)} \\
        & \xi_{j_k} = \xi_{j_{k-1}}
    \end{split}
    \]
    and it holds that
  \begin{eqnarray}\label{success_in}
   Z\left(\left[x_{j_{k-1}}+\frac{\alpha_{j_{k-1}}}{\delta^{s_{k-1}+1}} d_{k-1}\right]_{[l,u]};\eps\right) & \not < &  
   Z\left(\left[x_{j_{k-1}}+\frac{\alpha_{j_{k-1}}}{\delta^{s_{k-1}}} d_{k-1}\right]_{[l,u]};\eps\right)  \\
   && -\gamma\left(\left(\frac{\alpha_{j_{k-1}}}{\delta^{s_{k-1}+1}}\right)^2 - \left(\frac{\alpha_{j_{k-1}}}{\delta^{s_{k-1}}}\right)^2\right)\one.\nonumber
  \end{eqnarray}

  \item (Failure w.r.t. continuous variables) we have
      \begin{equation}\label{failure_in}
      Z([x_{j_{k-1}}+\alpha_{j_{k-1}} d_{k-1}]_{[l,u]};\eps) > Z(x_{\ell_j};\eps) -\gamma(\alpha_{j_{k-1}})^2\one,
      \end{equation}
     for at least a point $x_{\ell_j}\in \tilde L_{k-1}$, that is, $Z([x_{j_{k-1}}+\alpha_{j_{k-1}} d_{k-1}]_{[l,u]};\eps)$ is dominated by 
     $Z(x_{\ell_j};\eps) -\gamma(\alpha_{j_{k-1}})^2\one$. This implies that $[x_{j_{k-1}}+\alpha_{j_{k-1}}
    d_{k-1}]_{[l,u]}$ is a  ``bad'' point so that the step $\alpha_{j_{k-1}}$ associated to $x_{j_{k-1}}$ is shrunk by the constant factor $\theta$, i.e.
     \begin{subequations}\label{failure_step}
     \begin{eqnarray}
      x_{j_k} & = & x_{j_{k-1}} \\
      \alpha_{j_k} & = & \theta \alpha_{j_{k-1}}\\
      \alpha_{j_k}^{(d)} & = & \alpha_{j_{k-1}}^{(d)}\\
      \xi_{j_k} & = & \xi_{j_{k-1}};
     \end{eqnarray}
     \end{subequations}
     
    \item (Success w.r.t. discrete variables) we have
    \[
    \begin{split}
        & x_{j_k} = \left(x_{j_{k-1}} + \alpha_{j_{k},p_{k-1}}^{(d)} p_{k-1}\right);\\
        & \alpha_{j_k} = \alpha_{j_{k-1}} \\
        & \alpha_{j_k,p_{k-1}}^{(d)} \leq 2^{s_{k-1}}\alpha_{j_{k-1},p_{k-1}}^{(d)} \\
        & \xi_{j_k} = \xi_{j_{k-1}}
    \end{split}
    \]
    and it holds that
  \begin{eqnarray}\label{success_in_discr}
   Z\left(x_{j_{k-1}}+\beta p_{k-1};\eps\right) & \not < &  
   Z\left(x_{j_{k-1}} + \alpha_{j_{k},p_{k-1}}^{(d)} p_{k-1};\eps\right) 
   -\xi\one\\\nonumber
   \beta & \leq & { 2\alpha_{j_{k},p_{k-1}}^{(d)}}
  \end{eqnarray}

  \item (Failure w.r.t. discrete variables) we have
      \begin{equation}\label{failure_in_discr}
      Z([x_{j_{k-1}}+ p_{k-1}]_{[l,u]};\eps) > Z(x_{\ell_j};\eps) -\xi\one,
      \end{equation}
     for at least a point $x_{\ell_j}\in \tilde L_{k-1}$ and all $p_{k-1}\in D_{k-1}$. This implies 
     \begin{subequations}\label{failure_step_discr}
     \begin{eqnarray}
      x_{j_k} & = & x_{j_{k-1}} \\
      \alpha_{j_k} & = & \alpha_{j_{k-1}}\\
      \alpha_{j_k}^{(d)} & = & \alpha_{j_{k-1}}^{(d)} = 1\\
      \xi_{j_k} & = & \theta\xi_{j_{k-1}};
     \end{eqnarray}
     \end{subequations}
\end{enumerate}

\begin{lemma}[{\cite[Lemma 3.15]{liuzzi2016derivative}}]\label{proj_lemma}
Let $\alpha,\beta\in \Re$, with $\alpha, \beta \neq 0$, be such that $sign(\alpha)=sign(\beta)$ and let $x \in X$ and $p\in \Re^n$. Then, we have
\[
 \left[x+(\alpha+\beta)p \right]_{[l,u]}= \left[\left[x+\alpha p\right]_{[l,u]}+ \beta p\right]_{[l,u]}.
 \]
\end{lemma}

\begin{proposition}
The {\em Projected Expansion} is well-defined, that is, either the test at step \ref{LS:step3a} is satisfied or steps \ref{LS:step2}--\ref{LS:step12} are executed. In the latter case: 
\begin{itemize}
 \item[i)] the test at Step \ref{LS:step4} is satisfied a finite number of times, i.e. the procedure cannot infinitely cycle;
 \item[ii)] the procedure  updates $\tilde L$, i.e. the test at Step \ref{LS:step3} is satisfied at least once, i.e. {\tt Add\&Filter} is called at least once 
\end{itemize}
\end{proposition}

\par\smallskip\noindent
{\bf Proof}. If the test at step \ref{LS:step3a} is satisfied then, the procedure terminates. On the other hand, namely if the test at step \ref{LS:step3a} is not satisfied, the procedure executes steps \ref{LS:step2}--\ref{LS:step12} thus becoming equivalent to the {\em Projected expansion} in \cite{liuzzi2016derivative}. Hence, in the latter case, the proof follows from the proof of Proposition 3.16 of \cite{liuzzi2016derivative}. $\hfill\Box$

\par\medskip

In the following we prove an analogous result for the Discrete Search procedure.

\begin{proposition}
The {\em Discrete Search} is well-defined, that is, either the test at step \ref{DS:step4a} is satisfied or steps \ref{DS:step2}--\ref{DS:step15} are executed. In the latter case:
\begin{itemize}
 \item[i)] the test at Step \ref{DS:step4} is satisfied a finite number of times, i.e. the procedure cannot infinitely cycle;
 \item[ii)] the procedure updates $\tilde L$, i.e. either Step \ref{DS:step12} or \ref{DS:step19} are executed, that is {\tt Add\&Filter} is called at least once. 
\end{itemize}
\end{proposition}

\par\smallskip\noindent
{\bf Proof}. If the test at step \ref{DS:step4a} is satisfied, the procedure terminates. On the other hand, steps \ref{DS:step2}-\ref{DS:step15} get executed. 

\underline{Point i)}. We proceed by contradiction and assume that the test at Step \ref{DS:step4} is always satisfied, i.e. 
a monotonically increasing sequence of positive numbers $\{\beta^j\}$ exists. But this contradicts the fact that, by compactness of set $X\cap Z$, $\bar\alpha$ is finite. 
\par\smallskip

\underline{Point ii)}. If $\beta = \bar \alpha$, then the {\em Discrete Search} updates $\tilde L$ at step \ref{DS:step19} which concludes the proof.

Otherwise, let us suppose that $\beta < \bar\alpha$ at every iteration of the {\em Discrete Search}.. In this case,
let us consider a generic iteration. Then, we have
 \[
   \Big(\!Z([y+\alpha p]_{[l,u]};\eps) \not >  Z(x_j;\eps) -\xi\!\Big) \ \forall\ x_j\in \tilde L,
 \]
 either because $\alpha = \alpha_p^{(d)}$ or because the test at step \ref{DS:step4} is satisfied and $\alpha = \beta$. 
The above relation imply that, for all $x_j\in \tilde L$, and index $\ell_j\in\{1,\dots,q\}$ exists such that
\begin{equation}\label{Zellej}
   Z_{\ell_j}([y+\alpha p]_{[l,u]};\eps) \leq  Z_{\ell_j}(x_j;\eps) -\xi.
\end{equation}

Furthermore, let us assume that the test at Step \ref{DS:step3} is not satisfied, which means that
\[
 Z([y+\beta p]_{[l,u]};\eps) < Z([y+\alpha p]_{[l,u]};\eps) -\xi\one
\]
and, in particular, in view of (\ref{Zellej}), for all $x_j\in \tilde L$, and index $\ell_j\in\{1,\dots,q\}$ exists such that
\[
 Z_{\ell_j}([y+\beta p]_{[l,u]};\eps) < Z_{\ell_j}([y+\alpha p]_{[l,u]};\eps) - \gamma\xi \leq Z_{\ell_j}(x_j;\eps) -\xi.
\]
The above relation, recalling that $\beta <\bar\alpha$, implies that the test at Step \ref{DS:step4} is satisfied so that the procedure will perform a further iteration. 

Now, let us assume by contradiction that the test at Step \ref{DS:step3} is never satisfied. Then, the above reasoning would imply that the procedure infinitely cycles, 
which contradicts Point i). $\hfill\Box$

\begin{proposition}\label{mainconv}
Let Assumption \ref{assmfcq} hold. Let $\{L_k\}$  with $L_k = \{(x_j,\alpha_j),\ j=1,\dots,|L_k|\}$ be the sequence of sets of nondominated pairs produced by Algorithm DFMOINT.
Let $\{(x_{j_k},\alpha_{j_k})\}$ be a linked sequence and $\bar x$ be any limit point of $\{x_{j_k}\}$, i.e.,
\[
 \lim_{k\to\infty,k\in  K} x_{j_k} = \bar x,
\]
for a subset $K$ of indices. Then, an $\eps^\star > 0$ exists such that for all $\eps\in (0,\eps^\star]$,
if the subsequence $\{d_k\}_{k \in K}$ is dense in the unit sphere, $\bar x$
is a Pareto-Clarke KKT stationary point for Problem \eqref{prob} w.r.t. the continuous variables.
\end{proposition}

{\bf Proof}. If $\lim_{k\to\infty,k\in  K} x_{j_k} = \bar x$,
an index $\bar k$ exists such that, for all $k\geq\bar k$, $k\in K$, $(x_{j_k})_z = (\bar x)_z$.
Then, the proof follows from the proof of \cite[Proposition 3.19    ]{liuzzi2016derivative}. $\hfill\Box$

\par\medskip

\begin{proposition}
Let $\{L_k\}$ with $L_k=\{(x_j,\alpha_j,\alpha_j^{(d)},\xi_j), \ j=1,\dots, |L_k|\}$ be the
sequence of sets of nondominated points produced by DFMOINT. Every linked sequence $\{(x_{j_k},\alpha_{j_k},\alpha_{j_k}^{(d)},\xi_{j_k})\}$ is such that
\[
\lim_{k\to\infty}\xi_{j_k} = 0.
\]
\end{proposition}

{\bf Proof}. By the instruction of Algorithm DFNDFL, it follows that
$0< \xi_{j_{k+1}} \leq \xi_{j_k}$ for all $k$, meaning that the sequence $\{\xi_{j_k}\}$ is monotonically nonincreasing. Hence, $\{\xi_{j_k}\}$ converges to a limit $M \geq 0$. Suppose, by contradiction, that
$M > 0$. This implies that an index $\bar k > 0$ exists such that 
{
\begin{equation}\label{eq:xi_equal_M}
    \xi_{j_{k+1}} = \xi_{j_k} = M
\end{equation}
} for all $k\geq\bar k$.
Now, we preliminary show that, for every primitive direction $p$, a $\bar k$ exists such that, for all $k\geq \bar k$, $\alpha_{p,j_k}^{(d)} = 1$. We again proceed by contradiction and assume that and index set $K\subset\{0,1,2,\dots\}$ exists such that 
\[
  \alpha_{p,j_k}^{(d)} > 1,\quad\mbox{for all}\ k\in K.
\]
If this is the case, for $k\in K$, the test at step \ref{DS:step4a} of the Discrete Search procedure is not satisfied. This means that $\alpha > 0$ and, for all $x_j\in\tilde L_k$, $Z(x_{j_k}+\alpha p;\eps) \not > Z(x_j;\eps) -\xi_{j_k}\mathbf{1}$.  Hence, an index $\ell_k\in\{1,\dots,q\}$, exists such that
\[
  Z_{\ell_k}(x_{j_k}+\alpha p;\eps) \leq Z_{\ell_k}(x_{j_k};\eps) - \xi_{j_k} \leq  Z_{\ell_k}(x_{j_k};\eps) - M.
\]
Since, the number of objective functions is finite, it is possible to extract a further subset $\bar K\subseteq K$ such that, $\ell_k = \bar\ell$, for all $k\in\bar K$ and
\[
  Z_{\bar\ell}(x_{j_k}+\alpha p;\eps) \leq Z_{\bar\ell}(x_{j_k};\eps) - \xi_{j_k} \leq  Z_{\bar\ell}(x_{j_k};\eps) - M.
\]
The above relation is in contrast with continuity assumption of $Z_{\ell}$ and compactness of $X\cap \cal Z$.

Now, by compactness of $X\cap \cal Z$, the feasible primitive directions are finite. Hence, an index $\hat k$ exists such that
\[
   \alpha_{p,j_k}^{(d)} = 1,\quad\mbox{for all}\ p\in D,\ k\geq\hat k.
\]
Hence, for $k\geq \hat k$, the Discrete Search does not change the points within $L_{j_k}^{d}$, i.e. $L^+[x] = L_{j_k}^d[x]$. In this case, by the instruction of the algorithm DFMOINT, $\xi_{j_{k+1}} = \theta\xi_{j_k} = \theta M < \xi_{j_k}$ which is a contradiction with (\ref{eq:xi_equal_M}).
$\hfill\Box$

\par\medskip

For every linked sequence $\{(x_{j_k},\alpha_{j_k},\alpha_{j_k}^{(d)},\xi_{j_k})\}$, let us define the set 
\[
   H = \{k: \xi_{j_{k+1}} < \xi_{j_{k}}\}.
\]

\begin{proposition}
Let Assumption \ref{assmfcq} hold. Let $\{L_k\}$  with $L_k = \{(x_j,\alpha_j),\ j=1,\dots,|L_k|\}$ be the sequence of sets of nondominated pairs produced by Algorithm DFMOINT.
Let $\{(x_{j_k},\alpha_{j_k},\alpha_{j_k}^{(d)},\xi_{j_k})\}$ be a linked sequence and $\bar x$ be any limit point of $\{x_{j_k}\}_H$, i.e., a index $K\subseteq H$ exists such that
\[
 \lim_{k\to\infty,k\in  K} x_{j_k} = \bar x,
\]
Then, an $\eps^\star > 0$ exists such that for all $\eps\in (0,\eps^\star]$,
$\bar x$
is a Pareto-Clarke KKT stationary point for Problem \eqref{prob}.
\end{proposition}

{\bf Proof}. By the fact that $k\in K$, we know that an index $\bar k\in K$ exists such that

\begin{equation}\label{eq:discrchange}\begin{split}
    & (x_{j_k})_z = (x_{j_{k+1}})_z= \bar x_z,\quad\forall\ k\in K, k\geq\bar k,\\
    & (\alpha_{p}^{(d)})_{j_k} = 1,\quad\forall\ p\in D_k
 \end{split}
\end{equation}
Furthermore, since $\alpha_{j_k}^c \to 0$, we also have that 
\[
   \lim_{k\to\infty,k\in K}\alpha_{j_{k+1}}^c = 0.
\]
So that
\[
  \lim_{k\to\infty, k\in K}(x_{j_{k+1}})^c = (\bar x)^c,
\]
hence we get
\[
\lim_{k\to\infty,k\in K}x_{j_{k+1}} = \bar x.
\]

Let us consider any point $\hat x\in {\cal B}^z(\bar x)$. By the definition of discrete neighborhood ${\cal B}^z(\bar x)$, a direction $\hat d \in D^z(\bar x)$ exists such that
\begin{equation}\label{DFL_discr1-3}
\hat x = \bar x + \hat d.
\end{equation}
Then, {for all $k\in K$ sufficiently large}, {\eqref{eq:discrchange} and \eqref{DFL_discr1-3} imply}
\[
(x_{j_k} + \hat d)_z = (x_{j_{k+1}} + \hat d)_z = (\bar x  + \hat d)_z = (\hat x)_z.
\]
Hence, {for all $k\in K$ sufficiently large}, by the definition of discrete neighborhood we have $\hat d \in D^z(x_{j_{k+1}})$ and
$$
x_{j_{k}} + \hat d \in X\cap  \mathcal{Z}.
$$
Then, since $k \in K\subseteq H$, by the definition of $H$ and by the instructions of the Discrete search procedure, we have
\begin{equation*}
Z(x_{j_{k}} + \hat d;\eps) > Z(x_{j_k}) -\xi_{j_k}.
\end{equation*}
By taking the limit for $k\to\infty$, $k\in K$, we get
\[
Z(\hat x;\eps)\geq Z(\bar x),
\]
hence $\hat x$ cannot dominate $\bar x$ and the result follows by Proposition \ref{mainconv}. $\hfill\Box$

\section{Numerical Results}
In this section we report the results obtained in the numerical experimentation and comparison of the proposed algorithm with state-of-the-art methods in the literature on a set of multiobjective test problems.

\paragraph{Test problems description}\ \par
We used the 10 multiobjective problems proposed for the CEC 2009 special session and competition \cite{zhang2008multiobjective}. These problems have $n=10$ variables and two or three objective functions. Since the number of variables can be modified for all of the problems, we considered $n\in\{10,15,20,25,30\}$ thus obtaining a total of 50 problems.

In order to come up with mixed integer problems, given a continuous problem (among the 50 above defined) 
\[
  \begin{array}{l}
  \min\ \tilde f_1(\tilde x), \dots, \tilde f_q(\tilde x)\\
  s.t.\ \ell_i \leq \tilde x_i \leq u_i,\quad\forall\ i=1,\dots,n
  \end{array}
\]
we consider $n_c = \lfloor n/2\rfloor$ continuous variables and, consequently, $n_z = n-n_c$ integer variables. Hence, we define the following discretized problem
\begin{equation}\label{probnumbox}
  \begin{array}{l}
  \min\ f_1(x), \dots, f_q(x)\\
  s.t.\ \ell_i \leq x_i \leq u_i,\quad\forall\ i=1,\dots,n_c\\
  \phantom{s.t.}\ 0\leq x_i \leq 100,\quad\forall\ i=n_c+1,\dots,n
  \end{array}
\end{equation}
where 
\[
\begin{split}
  & \tilde x_i = x_i,\quad i=1,\dots,n_c,\\
  & \tilde x_i = \ell_i + \frac{u_i-\ell_i}{100}x_i,\quad i=n_c+1,\dots,n
\end{split}
\]
and $f_j(x) = \tilde f_j(\tilde x)$, $j=1,\dots,q$. 

As concerns constrained problems, given a bound constrained problem like problem (\ref{probnumbox}), we generate six constrained problems by adding the families of nonlinear constraints proposed in \cite{Karmitsa2007} and reported in table \ref{tab:constraints}.
\begin{table}[!ht]
\begin{center}\small
	\begin{tabular}{ll}
		$g_j(x) = (3-2x_{j+1})x_{j+1} - x_j -2x_{j+2} + 1 \leq 0,$ & for all \ $j \in \{1,{2,}\dots,n-2\}$,\\ \smallskip
		$g_j(x) = (3-2x_{j+1})x_{j+1} - x_j -2x_{j+2} + 2.5 \leq 0$, & for all \ $j \in \{1,{2,}\dots,n-2\}$,\\ \smallskip
		$ g_j( x) = x_j^2 + x_{j+1}^2 + x_jx_{j+1} -2x_j -2x_{j+1} +1 \leq 0$, & for all \ $j \in \{1,{2,}\dots,n-1\}$,\\ \smallskip
		$g_j(x) = x_j^2 + x_{j+1}^2 + x_jx_{j+1} -1 \leq 0$, & for all \ $j \in \{1,{2,}\dots,n-1\}$,\\ \smallskip
		$g_j( x) = (3-0.5x_{j+1})x_{j+1} - x_j -2x_{j+2} + 1 \leq 0$, & for all \ $j \in \{1,{2,}\dots,n-2\}$,\\ \smallskip
		$g_1( x) = \sum_{j=1}^{n-2}((3-0.5x_{j+1})x_{j+1} - x_j -2x_{j+2} + 1) \leq 0$. & $\phantom{j=\{1,{2,}\dots,n-2\}}$.
	\end{tabular}
\end{center}
\caption{Families of constraints from \cite{Karmitsa2007}.}\label{tab:constraints}
\end{table}
Hence, we come up with 300 nonlinearly constrained multiobjective problems with number of variables $n\in\{10,15,20,25,30\}$, number of objective functions $q\in\{2,3\}$ and number of constraints $m\in\{1,8,9,13,14,18,19,23,24,28,29\}$. 

\paragraph{Methods}\ \par
The proposed method DFMOINT implemented in python 3.9 (along with the test problems) is freely and publicly available on the Derivative-Free Library DFL at the URL 

\centerline{\tt http://www.iasi.cnr.it/$\sim$liuzzi/DFL/} 
In the experimentation, we compared DFMOINT with two representative methods of, respectively, deterministic and probabilistic algorithms for mixed-integer multiobjective black-box optimization problems. Namely, we considered: 
\begin{itemize}
    \item the  {\tt BiMADS} algorithm for bi-objective problems implemented in the NOMAD (version 3.9.1) software package \cite{le2011algorithm,nomadwebsite,audet2008multiobjective};
    \item the NSGA-II algorithm \cite{deb2002fast} implemented in the {\tt pygmo} python library.
\end{itemize}
All of the methods were run by using their respective default parameters and allowing a maximum number of 20,000 function evaluations. In particular, NSGA-II was run using a population size of 100 individuals and allowing 200 generations so that NSGA-II uses exactly 20,000 function evaluations.

Since NSGA-II is a probabilistic algorithm, its outcomes vary if it is run using different seeds for the random sequence. Hence, we run NSGA-II 10 times with ten different random seed. Then, on the basis of the results obtained by these ten runs, we extracted the best, median and worst (ideal) algorithms for every single problem. The best run simply consists of the
run that has the higher percentage of nondominated solutions when compared to the
remaining ones (considering as a reference Pareto front the one obtained from the 10
runs performed). In a similar way, the median run is the one with a percentage of nondominated points that is the median one. Finally, the worst run is selected as the one with the lowest percentage of nondominated points.

\paragraph{Metrics and performance profiles}\ \par
The results of the various comparisons are reported in terms of the purity \cite{bandyopadhyay2004multiobjective},  $\Gamma$ and $\Delta$ spread metrics
(as defined in \cite{custodio2011direct}) by using
performance profiles \cite{dolan2002benchmarking}.

\begin{figure*}[t!]
    \centering
    \begin{subfigure}[t]{0.4\textwidth}
    \centering
    \includegraphics[width=0.98\textwidth]{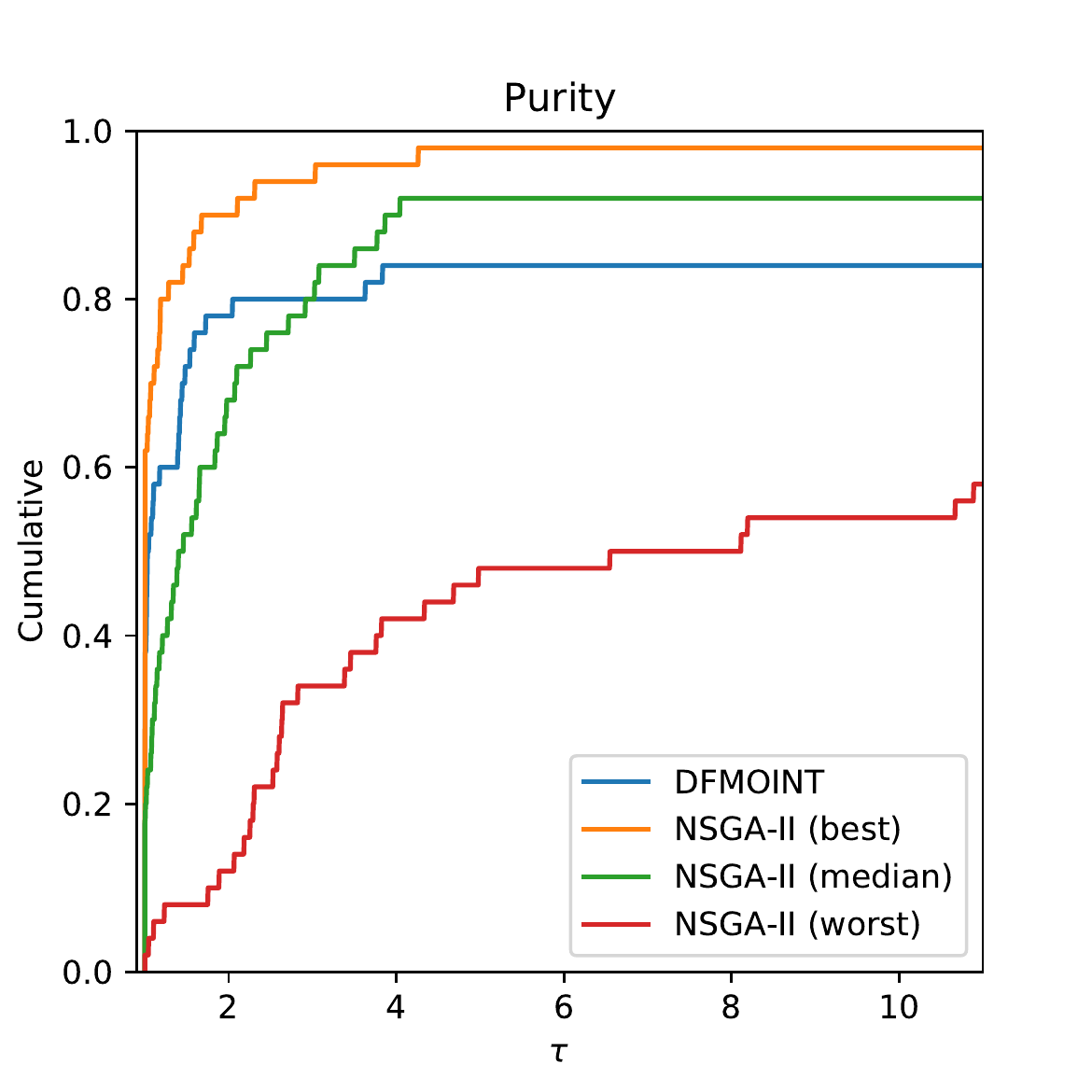}
    \caption{Purity performance profiles}
    \label{fig:purity}
    \end{subfigure}
    \begin{subfigure}[t]{0.4\textwidth}
    \centering
    \includegraphics[width=0.98\textwidth]{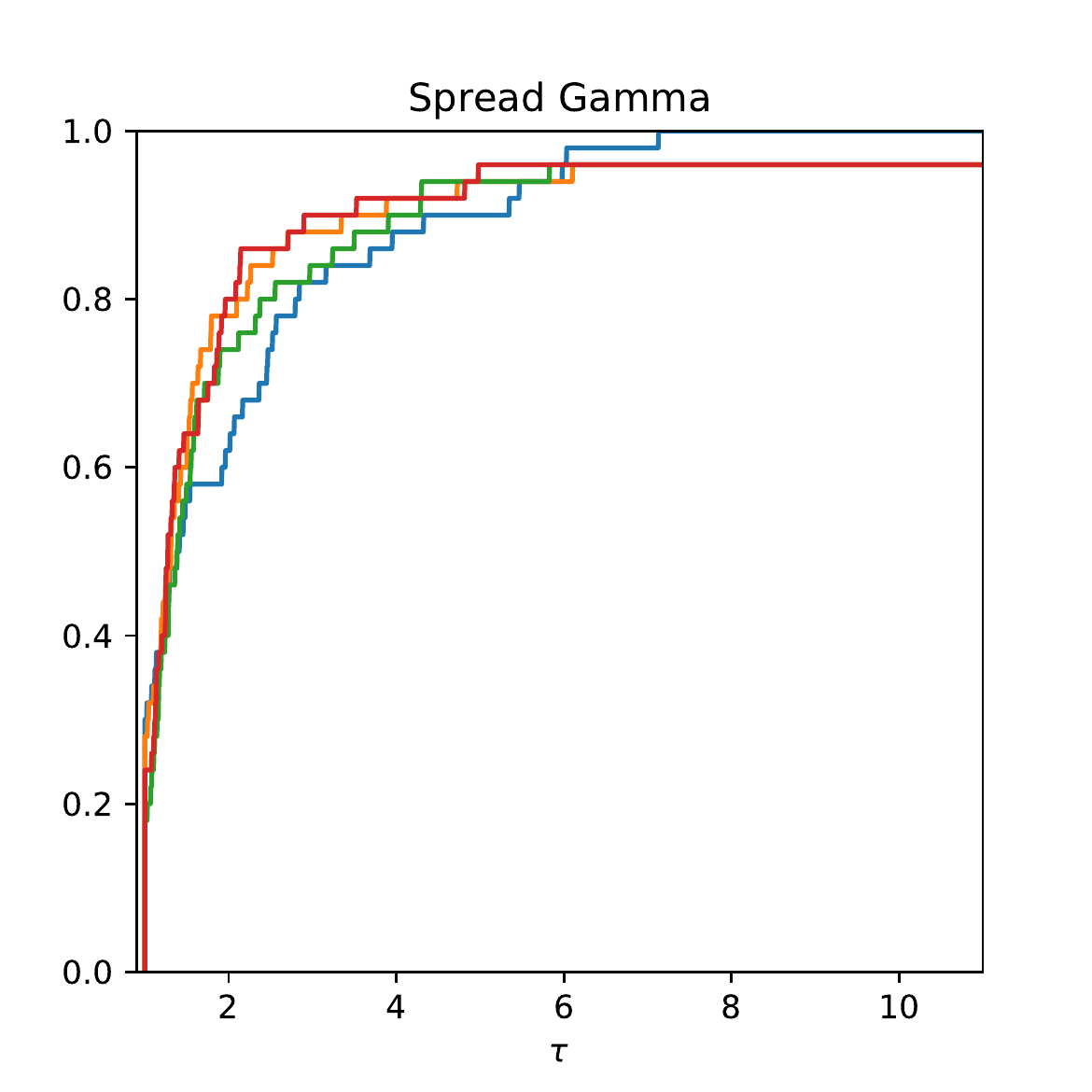}
    \caption{Spread $\Gamma$ performance profiles}
    \label{fig:spreadg}
    \end{subfigure}
    \begin{subfigure}[t]{0.4\textwidth}
    \centering
    \includegraphics[width=0.98\textwidth]{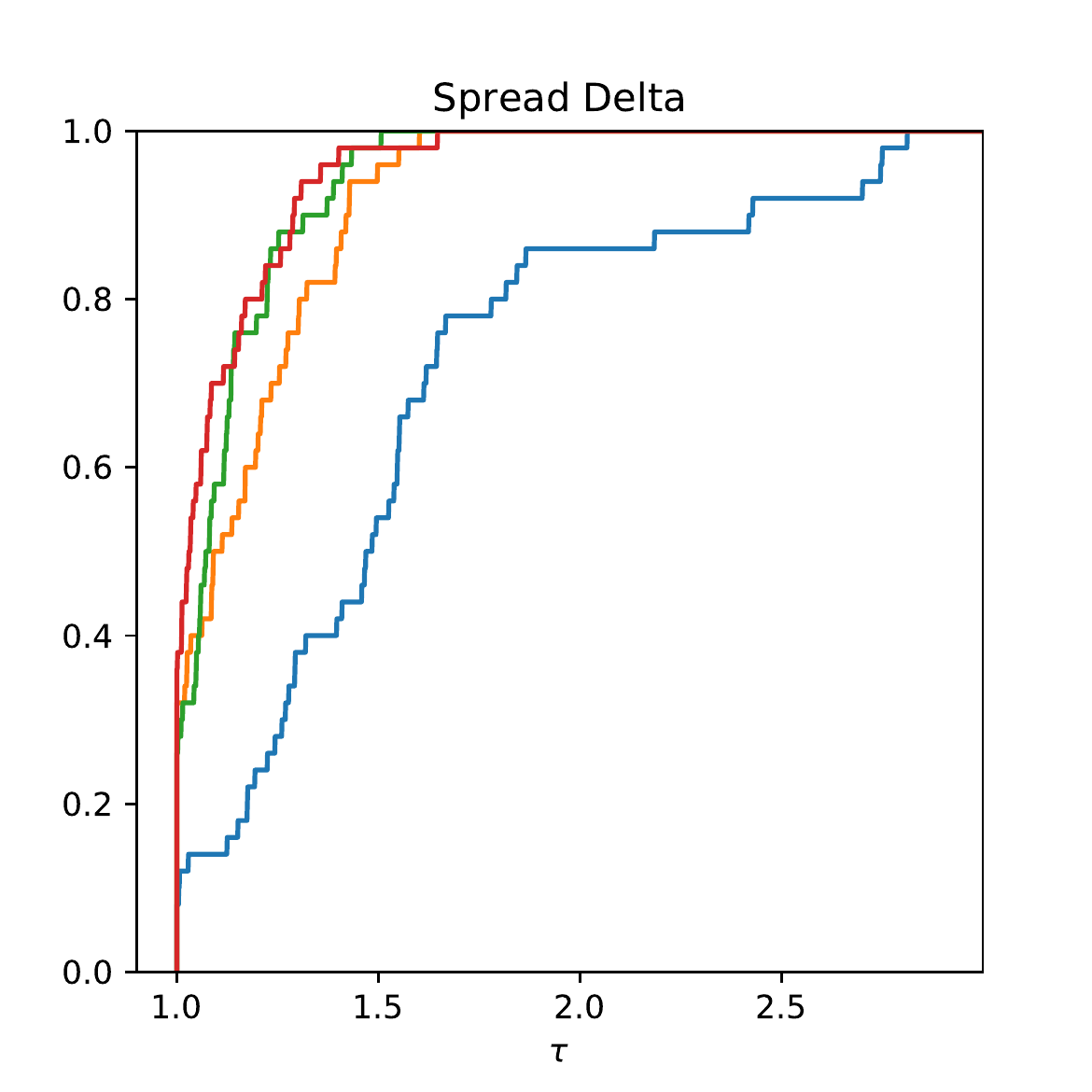}
    \caption{Spread $\Delta$ performance profiles}
    \label{fig:spreadd}
    \end{subfigure}
\caption{Comparison of DFMOINT with NSGA-II best, worst and median on the 50 bound constrained problems in terms of purity, spread $\Gamma$ and $\Delta$ performance profiles.}
    \label{fig:0}
\end{figure*}

\begin{figure*}[t!]
    \centering
    \begin{subfigure}[t]{0.4\textwidth}
    \centering
    \includegraphics[width=0.98\textwidth]{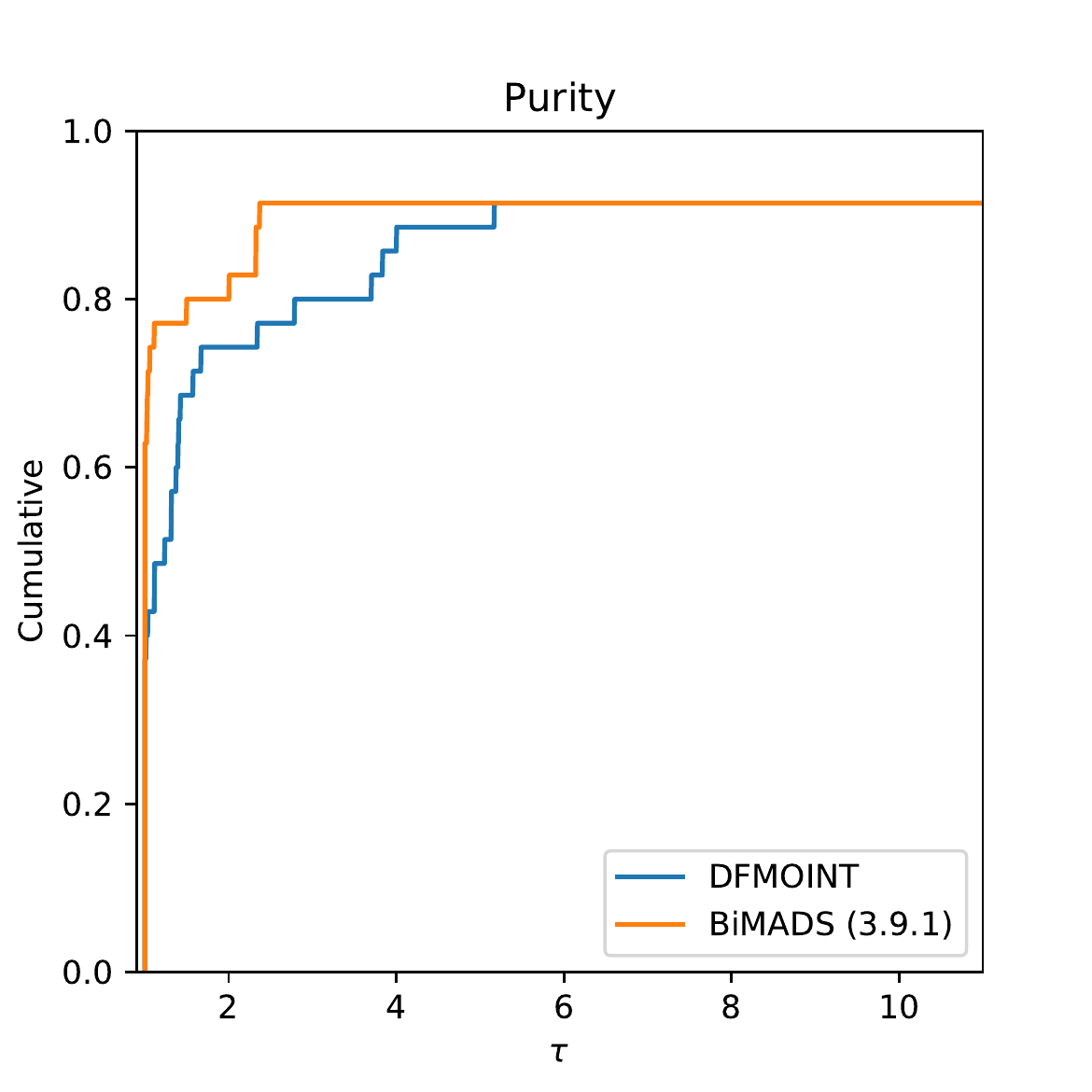}
    \caption{Purity performance profiles}
    \label{fig:purity1}
    \end{subfigure}
    \begin{subfigure}[t]{0.4\textwidth}
    \centering
    \includegraphics[width=0.98\textwidth]{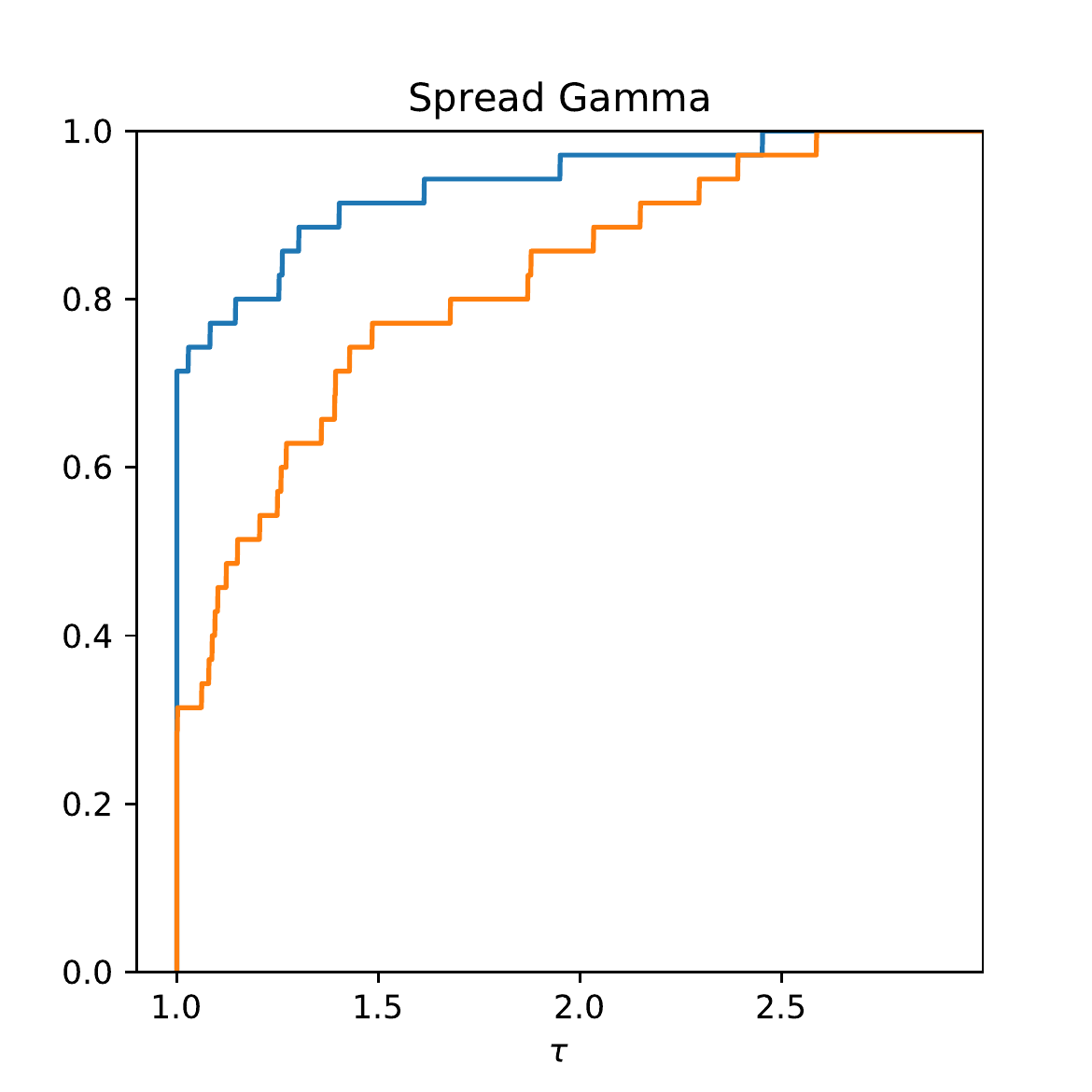}
    \caption{Spread $\Gamma$ performance profiles}
    \label{fig:spreadg1}
    \end{subfigure}
    \begin{subfigure}[t]{0.4\textwidth}
    \centering
    \includegraphics[width=0.98\textwidth]{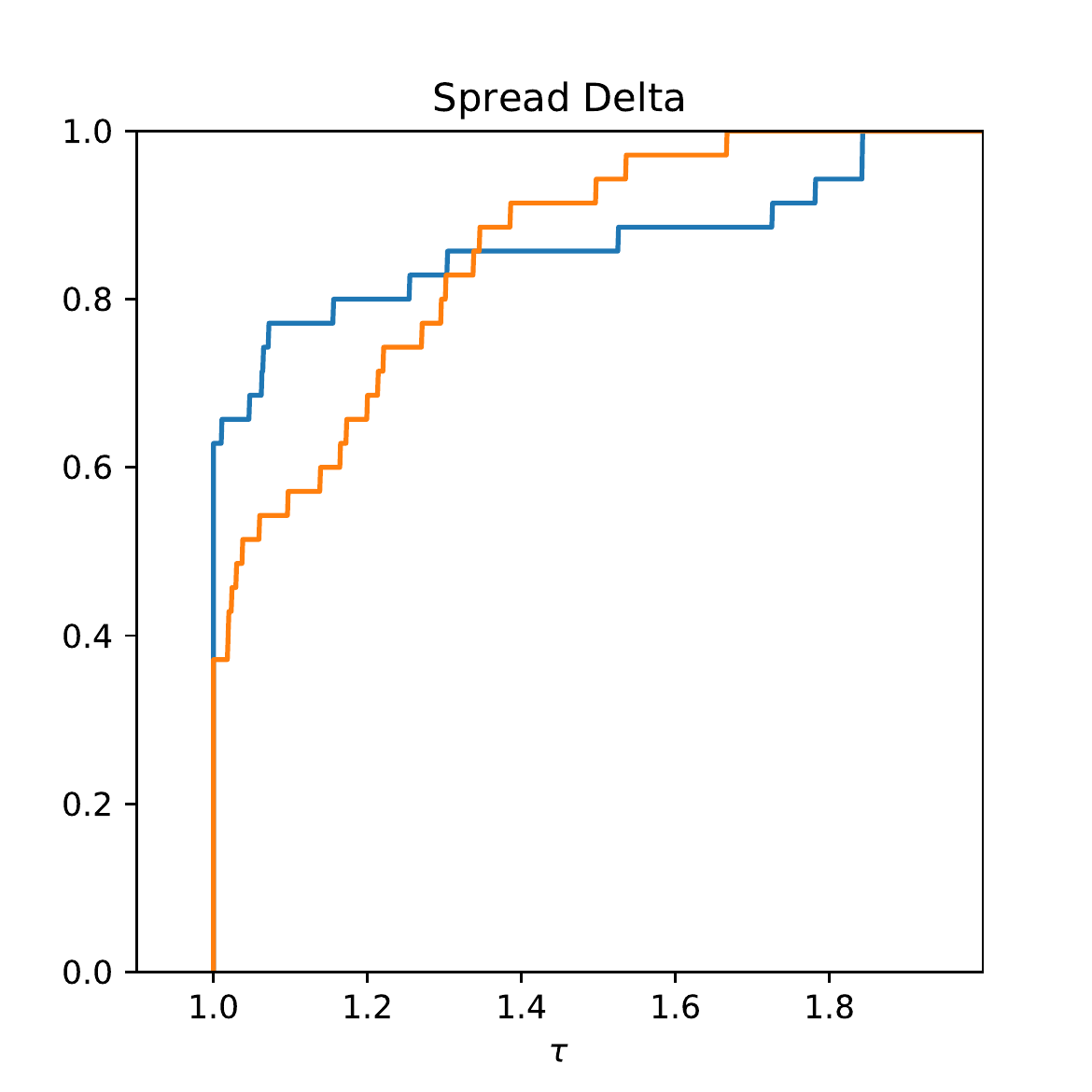}
    \caption{Spread $\Delta$ performance profiles}
    \label{fig:spreadd1}
    \end{subfigure}
\caption{Comparison of DFMOINT with BiMADS (implemented in NOMAD 3.9.1) on the 35 bound constrained problems with $q=2$ objective functions in terms of purity, spread $\Gamma$ and $\Delta$ performance profiles.}
    \label{fig:1}
\end{figure*}

\paragraph{Results on bound constrained problems}
In figure \ref{fig:0} we report the comparison of the proposed algorithm DFMOINT with the three versions of NSGA-II, namely the best, worst and median. The results are reported in terms of performance profiles. As we can see, DFMOINT compares favorably with NSGA-II (median) in terms of purity. In particular, DFMOINT is more efficient than NSGA-II (median) but slightly less robust. As regards, the spread $\Gamma$ metric, all of the algorithms have very similar performances whereas, for the spread $\Delta$ metric DFMOINT behaves poorly with respect to NSGA-II. If we consider that we are comparing our algorithm against NSGA-II which is designed to perform a global exploration of the search space, the comparison is not that bad. In fact, the following observations can be carried out:
\begin{itemize}
    \item[(i)] in terms of purity, DFMOINT performances are closer to the best version of NSGA-II rather than to the worst version; furthermore, DFMOINT is competitive with the median version of NSGA-II;
    \item[(ii)] the behaviour of DFMOINT with respect to the spread metrics was largely expected given the local optimization nature of our proposed algorithm and again considering the definition of NSGA-II.
\end{itemize}
Concerning the comparison with BiMADS, we have to first remark that BiMADS can only solve problems with only two objective functions, whereas DFMOINT is able to manage an arbitrary number of them. Hence, the comparison has been carried out on the subset of problems with $q=2$ objective functions, i.e. 35 problems out of the 50 bound constrained problems.
The results are reported in figure \ref{fig:1}. DFMOINT compares favorably with BiMADS in terms of spread $\Delta$ and is the clear winner w.r.t. the spread $\Gamma$. For what concerns purity, it can be noted that BiMADS is more efficient than DFMOINT (on the subset of two-objective problems) but the two methods have the same robustness. It must also be said that BiMADS has been run using default values for the parameters. In particular, the scalarized problems are solved with MADS using Nelder-Meade heuristic in the initial phase of the optimization process. 

\paragraph{Results on constrained problems}
For the constrained problems set, it must be said that the adding of constraints reported in table \ref{tab:constraints} can conflict with the bound constraints on the variables.  Indeed, we noticed that for some problems no solver is able to find feasible non-dominated points. Hence, in the comparison we keep only those problems where at least one solver is able to find feasible non-dominated points.

\begin{figure*}[t!]
    \centering
    \begin{subfigure}[t]{0.4\textwidth}
    \centering
    \includegraphics[width=0.98\textwidth]{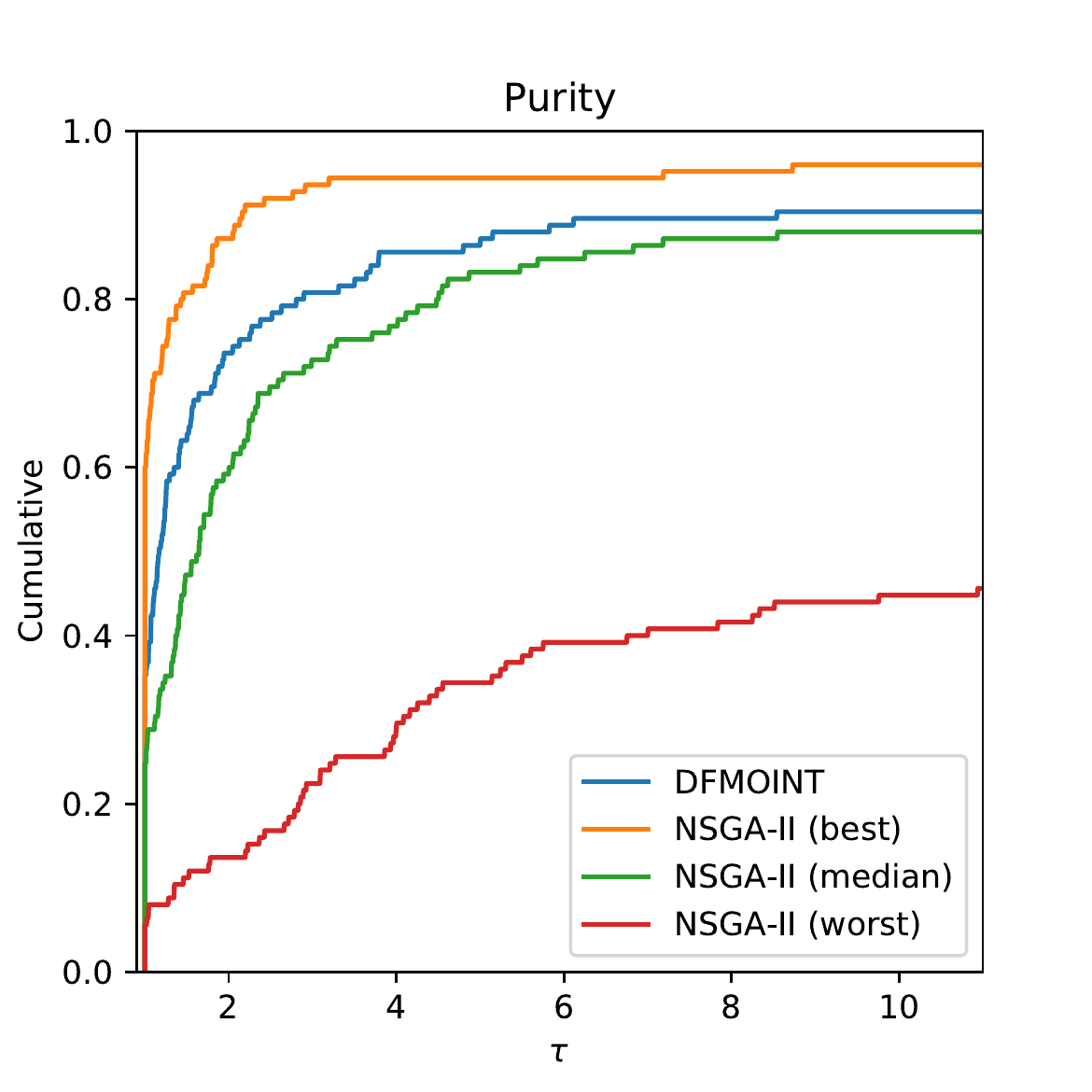}
    \caption{Purity performance profiles}
    \label{fig:purity2}
    \end{subfigure}
    \begin{subfigure}[t]{0.4\textwidth}
    \centering
    \includegraphics[width=0.98\textwidth]{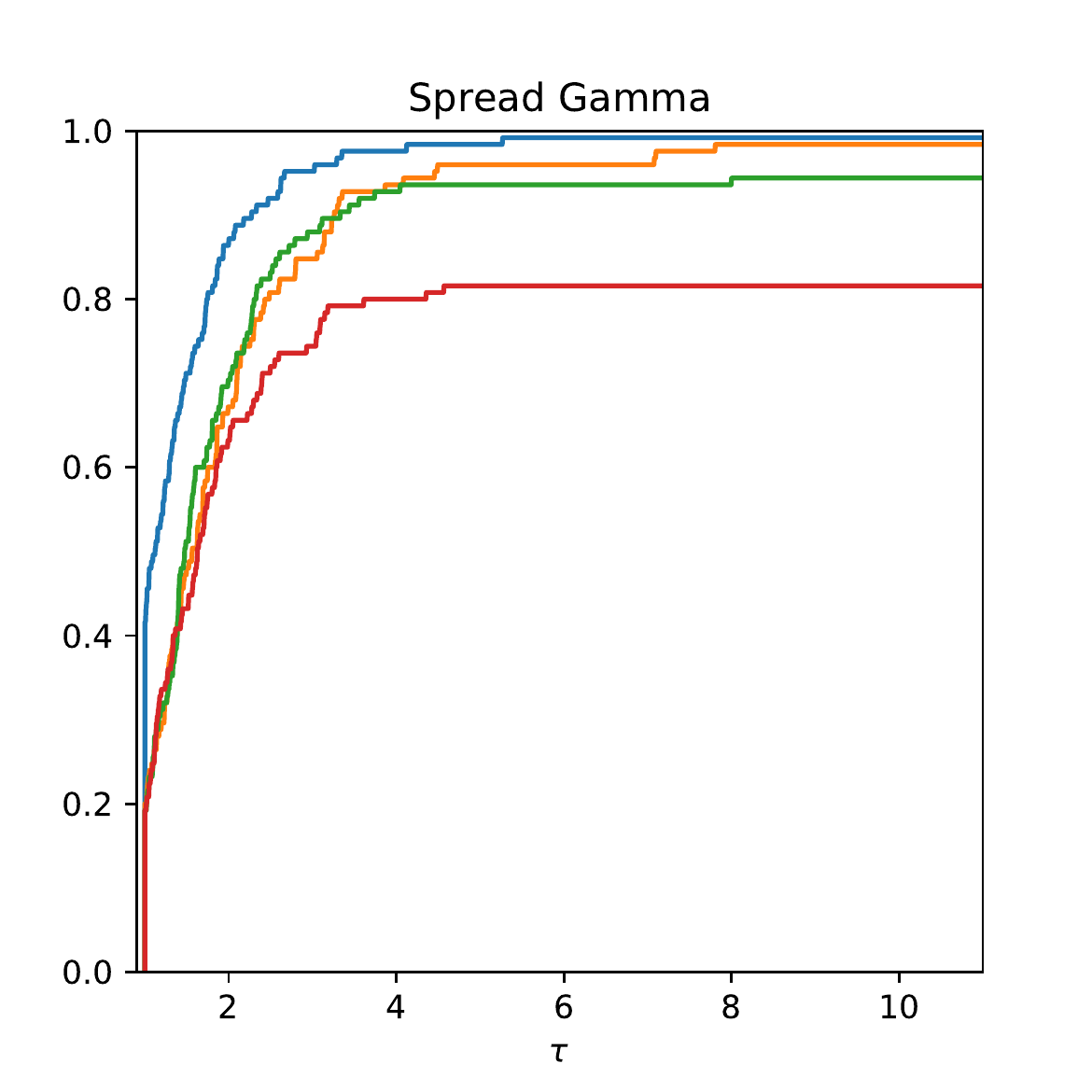}
    \caption{Spread $\Gamma$ performance profiles}
    \label{fig:spreadg2}
    \end{subfigure}
    \begin{subfigure}[t]{0.4\textwidth}
    \centering
    \includegraphics[width=0.98\textwidth]{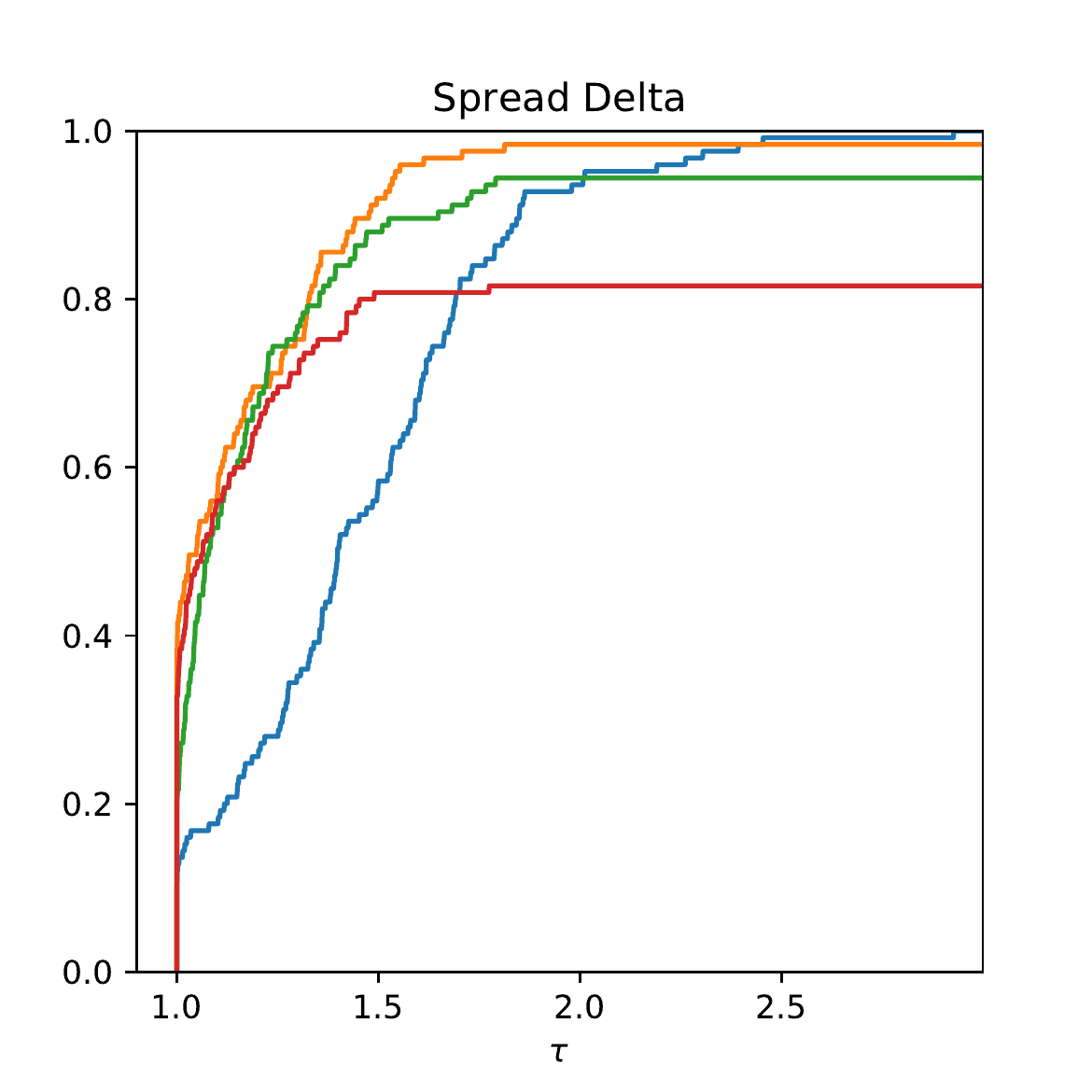}
    \caption{Spread $\Delta$ performance profiles}
    \label{fig:spreadd2}
    \end{subfigure}
\caption{Comparison of DFMOINT with NSGA-II best, worst and median on the subset of 300 constrained problems where at least one solver finds feasible nondominated points in terms of purity, spread $\Gamma$ and $\Delta$ performance profiles.}
    \label{fig:2}
\end{figure*}

\begin{figure*}[t!]
    \centering
    \begin{subfigure}[t]{0.4\textwidth}
    \centering
    \includegraphics[width=0.98\textwidth]{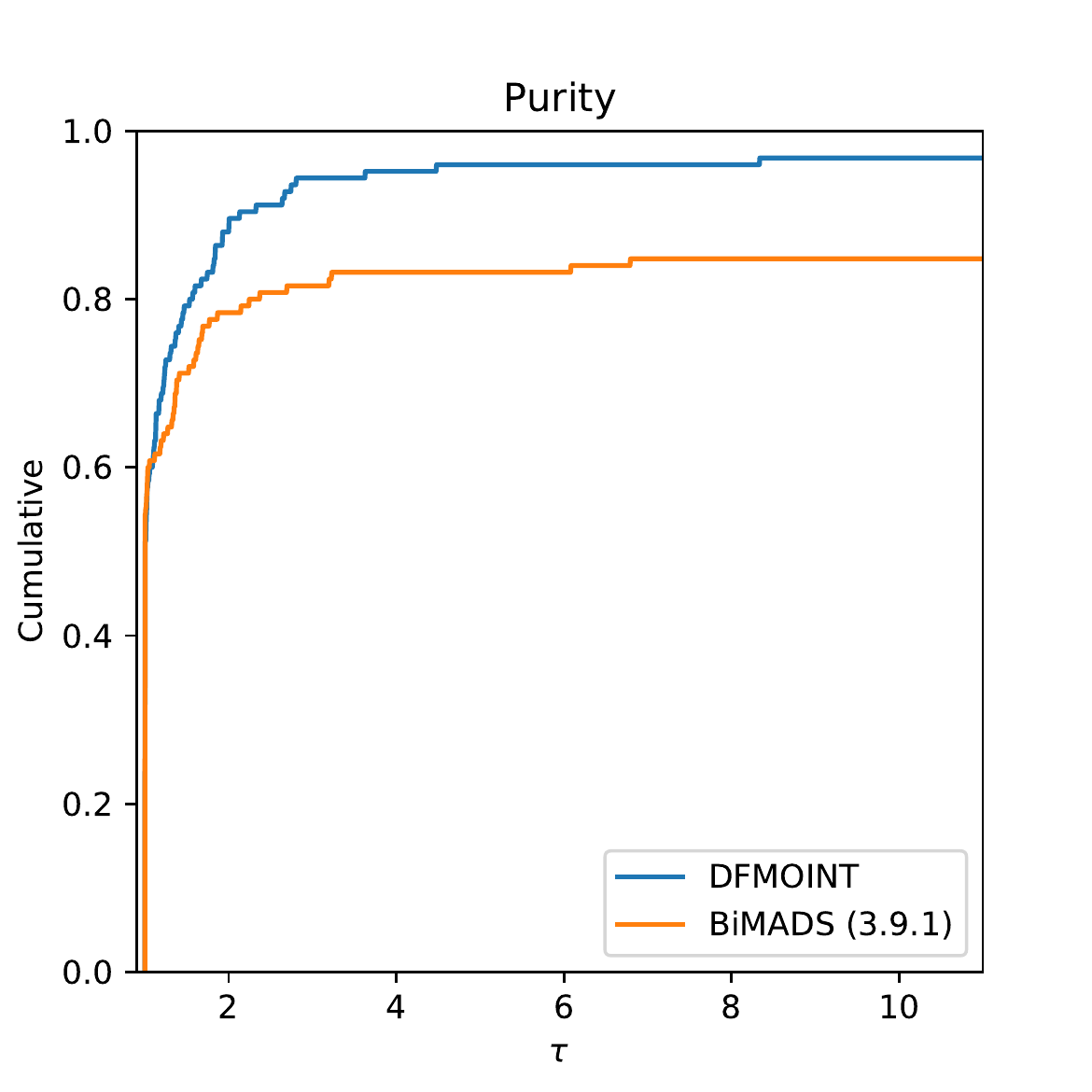}
    \caption{Purity performance profiles}
    \label{fig:purity3}
    \end{subfigure}
    \begin{subfigure}[t]{0.4\textwidth}
    \centering
    \includegraphics[width=0.98\textwidth]{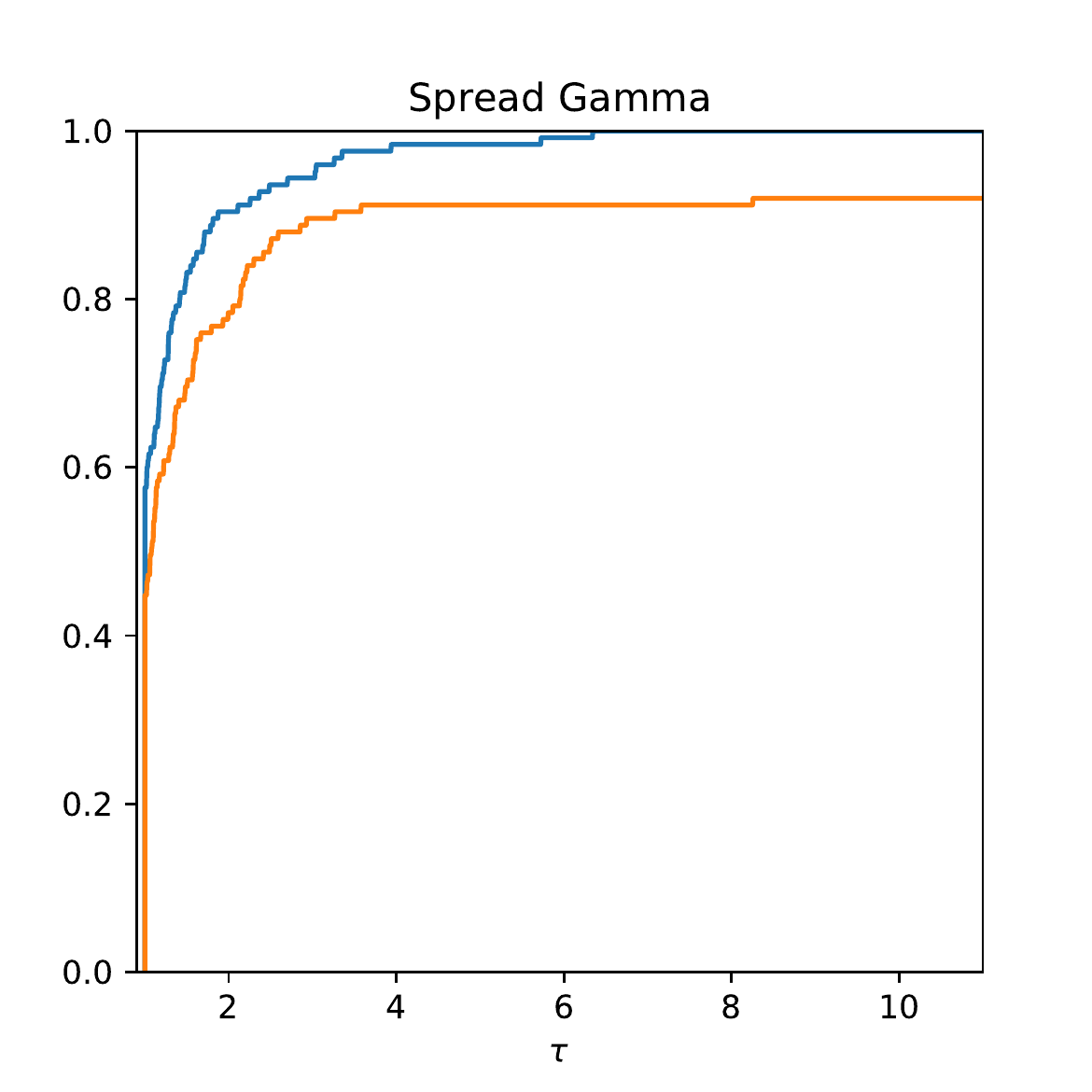}
    \caption{Spread $\Gamma$ performance profiles}
    \label{fig:spreadg3}
    \end{subfigure}
    \begin{subfigure}[t]{0.4\textwidth}
    \centering
    \includegraphics[width=0.98\textwidth]{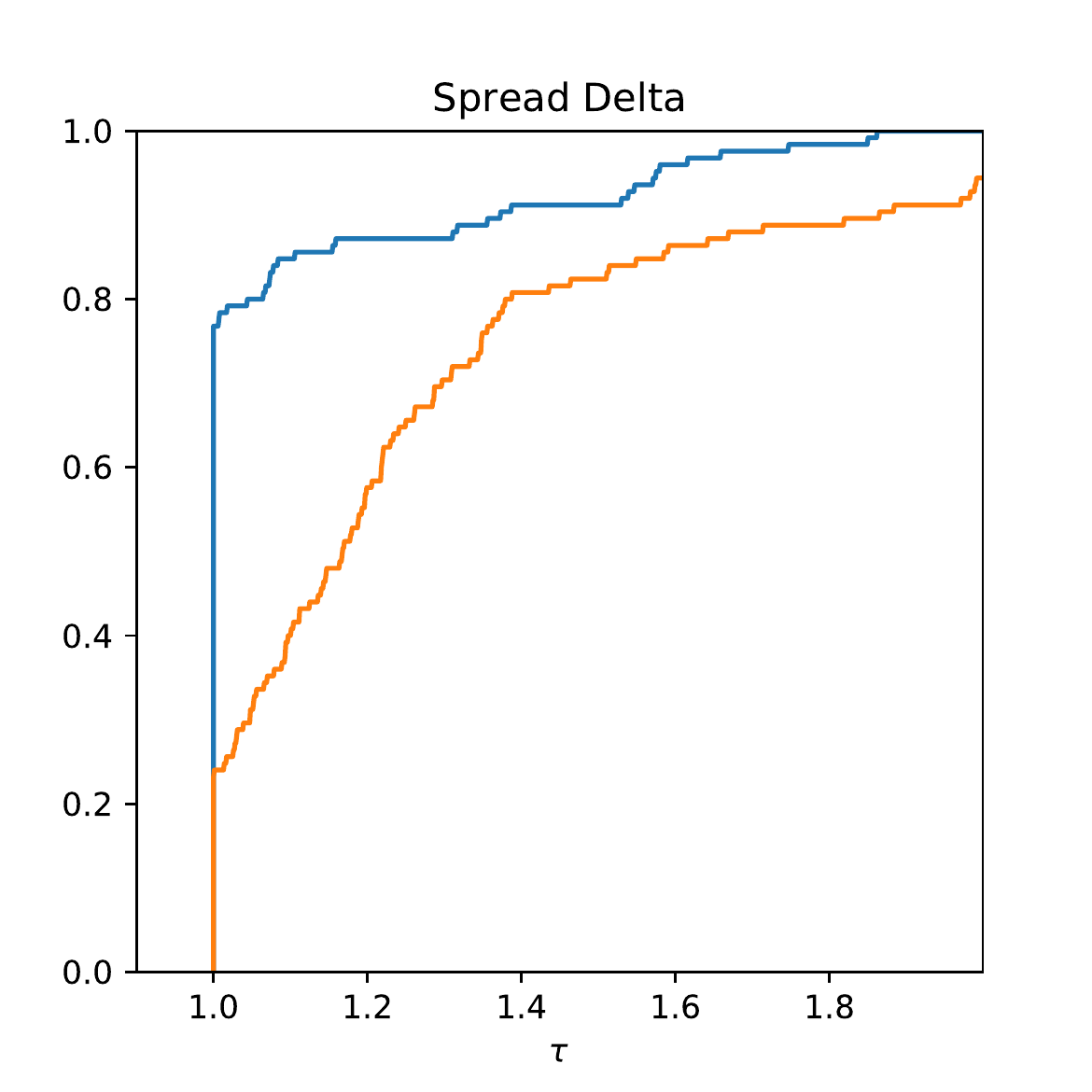}
    \caption{Spread $\Delta$ performance profiles}
    \label{fig:spreadd3}
    \end{subfigure}
\caption{Comparison of DFMOINT with BiMADS (implemented in NOMAD 3.9.1) on the subset of 210 bi-objective constrained problems where at least one solver finds feasible nondominated points in terms of purity, spread $\Gamma$ and $\Delta$ performance profiles.}
    \label{fig:3}
\end{figure*}

As concerns the behavior of DFMOINT on the nonlinearly contrained problems, figure \ref{fig:2} reports the comparison of DFMOINT with respect to the three versions of NSGA-II. We can make the following observations:
\begin{itemize}
    \item[(i)] in terms of purity, DFMOINT is better than NSGA-II median and worst and it is slightly dominated by the best version of NSGA-II;
    \item[(ii)] in terms of spread $\Gamma$, DFMOINT is the clear winner with respect to all of the versions of NSGA-II;
    \item[(iii)] in terms of spread $\Delta$, NSGA-II (each version) is more efficient than DFMOINT but less robust.
\end{itemize}
More or less the same observations can be made with regards to the comparison between DFMOINT and BiMADS on the subset of the 210 bi-objective constrained problems where at least one solver finds feasible non-dominated points. The profiles for such a comparison are reported in figure \ref{fig:3}. As we can see, DFMOINT is better than BiMADS with respect to all of the three metrics considered, namely purity, spread $\Gamma$ and $\Delta$.

\section{Conclusions}
In this paper we proposed a new derivative-free algorithm for the solution of mixed-integer black-box constrained multiobjective optimization problems. In particular, considering the presence of integer or discrete variables and the possible nonsmoothness of the functions defining the problem, we give local/global Pareto optimality definitions and stationarity conditions. 

Then, through the use of an exact penalty approach we show the equivalence between the original nonlinearly constrained multiobjective problem and a penalized problem where only bound constraints are present.

For the proposed algorithm, we carry out a thorough convergence analysis and we manage to prove convergence toward stationary points of the problem. 

The proposed algorithm (DFMOINT) has been tested on a large number of bound constrained problems and nonlinearly constrained problems and its performances have been compared with other state-of-the-art solvers for mixed-integer derivative-free multiobjective optimization. In particular, we considered in the experimentation two well-known solvers, namely BiMADS (implemented in the software package NOMAD v3.9.1) and NSGA-II (implemented in the pygmo python library). The results and comparison show the efficiency of DFMOINT with respect to both probabilistic and deterministic approaches.


\end{document}